\documentclass[12pt,oneside,reqno]{amsart}
\usepackage{amsfonts,amsmath,amscd,amssymb,paralist}
\usepackage{mathrsfs}
\usepackage{amsxtra}
\usepackage[mathscr]{eucal}
\usepackage{graphics}
\usepackage{latexsym,graphicx,verbatim,nameref}
\usepackage[numbers,sort&compress]{natbib}
\usepackage[all,cmtip]{xy}
\usepackage{hyperref}
\usepackage{fancyhdr}
\usepackage{color}
\usepackage{tikz}
\usepackage{tikz-cd}
\usepackage{graphicx}
\usetikzlibrary{arrows,scopes,matrix}
\numberwithin{equation}{section}

\makeatletter
\newtheorem*{rep@theorem}{\rep@title}
\newcommand{\newreptheorem}[2]{%
\newenvironment{rep#1}[1]{%
 \def\rep@title{#2 \ref{##1}}%
 \begin{rep@theorem}}%
 {\end{rep@theorem}}}
\makeatother

\newreptheorem{theorem}{Theorem}
\newreptheorem{corollary}{Corollary}
\newreptheorem{lemma}{Lemma}
\newreptheorem{conjecture}{Conjecture}

\theoremstyle{definition}

\newtheorem*{theorem*}{Theorem}
\newtheorem*{proposition*}{Proposition}

\newtheorem*{claim*}{Claim}

\newtheorem*{conjecture*}{Conjecture}
\newtheorem*{observation*}{Observation}
\newtheorem*{question*}{Question}

\begin{document}
\title{Compact Group Homeomorphisms Preserving The Haar Measure}
\author{Gang Liu}
\address{College of Mathematics and Statistics, Chongqing University}
\email{20181890@cqu.edu.cn}
\date{\today}
\begin{abstract}
This paper studies the measure-preserving homeomorphisms on compact groups and proposes new methods for constructing measure-preserving homeomorphisms on direct products of compact groups and non-commutative compact groups.

On the direct product of compact groups, we construct measure-preserving homeomorphisms using the method of integration. In particular, by applying this method to the \(n\)-dimensional torus \({\mathbb{T}}^{n}\), we can construct many new examples of measure-preserving homeomorphisms. We completely characterize the measure-preserving homeomorphisms on the two-dimensional torus where one coordinate is a translation depending on the other coordinate, and generalize this result to the \(n\)-dimensional torus.

For non-commutative compact groups, we generalize the concept of the normalizer subgroup \(N\left( H\right)\) of the subgroup \(H\) to the normalizer subset \({E}_{K}( P)\) from the subset \(K\) to the subset \(P\) of the group of measure-preserving homeomorphisms.  We prove that if \(\mu\) is the unique \(K\)-invariant measure, then the elements in \({E}_{K}\left( P\right)\) also preserve \(\mu\). In some non-commutative compact groups the normalizer subset \({E}_{G}\left( {\mathrm{{AF}}\left( G\right) }\right)\) can give non-affine homeomorphisms that preserve the Haar measure. Finally, we prove that when \(G\) is a finite cyclic group and a \(n\)-dimensional torus, then \(\mathrm{{AF}}\left( G\right)= N\left( G\right)  = {E}_{G}\left( {\mathrm{{AF}}\left( G\right) }\right)\).
\end{abstract}

\maketitle
\tableofcontents

\section{Introduction}\

Homeomorphisms preserving the Haar measure have been widely studied and applied in multiple fields. They have important applications in group theory, harmonic analysis, algebraic topology, geometry, and they occupy an important position especially in the research of dynamical systems.

Entropy is one of the most important concepts in dynamical systems. Generally, entropy can be divided into measure entropy and topological entropy. According to the variational principle, in a topological dynamical system, the topological entropy is equal to the supremum of the measure entropies of all invariant measures. The measure that makes the measure entropy equal to the topological entropy (i.e., the measure of maximal entropy) is a particularly concerned issue in dynamical systems. And the homeomorphism that preserves the Haar measure plays an extremely important role in the measure of maximal entropy. K. R. Berg\cite{1} proved that when the Haar measure is adopted in a topological dynamical system, the measure entropy reaches the maximum. Moreover, if the relevant transformation is ergodic and the entropy is finite, then the Haar measure is the only measure with this property. In the article of C. Deninger\cite{2}, the significance of the Haar measure in the dynamical system induced by the action of an amenable group is mentioned. He proved that for a compact metric abelian group \(G\) with a normalized Haar measure \({m}_{G}\), if a discrete amenable group \(H\) acts on \(G\) by automorphisms, then the topological entropy of the dynamical system (G, H) is still equal to the measure entropy of \(\left( {G,H,{m}_{G}}\right)\). From C. Deninger's proof, it can be seen that as long as a discrete amenable group \(H\) acts on a compact metric abelian group \(G\) with a normalized Haar measure \({m}_{G}\) by a homeomorphism that preserves the Haar measure, the measure entropy of \(\left( {G,H,{m}_{G}}\right)\) is still equal to the topological entropy of the dynamical system (G, H). Further, if \(T\) is a surjective homomorphism of the compact metric group \(G\) , it can still be proved that the Haar measure is an invariant measure of \(T\), and the entropy of the system reaches the maximum value when the Haar measure is taken, which is equal to the topological entropy of the topological dynamical system.

Some progress has also been made in the research on homeomorphisms that preserve the Haar measure other than entropy. As early as 1954, R. Arens\cite{3} conducted research on this. We know that homeomorphisms that preserve the Haar measure on the real line are actually those that preserve the metric. It is easy to know that such homeomorphisms have the form of \(f\left( x\right)  = {ax} + b\), that is, \(f\) is affine (a combination of translation and automorphism). In his work, R. Arens used the concept of quasi-invariant measure and the connectivity to extend the result that homeomorphisms preserving the Haar measure on the real line are affine to locally compact one-dimensional connected abelian groups. This result still has practical value in the recent number-theoretic research of the \(p\)-adic number field. Not all homeomorphisms that preserve the Haar measure have such a simple structure. The homeomorphism groups that preserve the Haar measure of many compact groups have very complex group structures. For these complex cases, we can only turn to studying their properties. In 2007, A. S. Kechris and C. Rosendal\cite{4} found that when studying the density and conjugacy of the automorphism groups of some countable structures (the automorphisms of compact groups preserve the Haar measure), the homeomorphism group \(\operatorname{Homeo}\left( {{2}^{\mathbb{N}},{m}_{{2}^{\mathbb{N}}}}\right)\) that preserves the Haar measure on the set \({2}^{\mathbb{N}}\) with the uniform convergence topology is a Polish group (a separable metric group) and it has a typical conjugacy class. In 2009, C. Rosendal\cite{5} proved on the basis of this research that the typical elements in \(\operatorname{Homeo}\left( {{2}^{\mathbb{N}},{m}_{{2}^{\mathbb{N}}}}\right)\) are conjugate to their non-zero powers, and further obtained that the typical elements in the isometry group of any rational Urysohn metric group (isometric mappings preserve the Haar measure) are also conjugate to their non-zero powers. Some progress has also been made in the research on homeomorphisms that preserve the Haar measure on the torus. In 2016, M. Andersson\cite{6} used the homotopy between homeomorphisms and linear mappings on the torus to prove that a homeomorphism that preserves the Haar measure and has a mapping degree greater than or equal to 2 on the two-dimensional torus is non-transitive if and only if the eigenvalues of its homotopic linear mapping have \(\pm  1\), that is, the eigenvalues of the induced matrix of the linear mapping have \(\pm  1\).

As an invariant measure for homeomorphisms that preserve the Haar measure, whether the Haar measure is unique, that is, whether the invariant measure of homeomorphisms that preserve the Haar measure is only the Haar measure in the sense of equivalence, is also one of the research directions. It is easy to see from Arens' results that the maps on the unit circle that preserve the Haar measure are only the compositions of translations and conjugations. However, there are many measure-preserving homeomorphisms for measures equivalent to the Haar measure. Whether the invariant measures of these homeomorphisms are equivalent to the Haar measure is an important problem. As early as 1981, M. Herman\cite{7} proved that for a \({C}^{2}\) diffeomorphism \(f\) on the circle with an irrational rotation number \(\alpha\), if \(F\) approximates an irrational rotation according to \(\alpha\), then \(f\) has an invariant measure equivalent to the Haar measure. In 1989, Y. Katznelson and D. Ornstein\cite{8} improved this result by removing the approximation condition. In 2004, I. Liousse \cite{9} further proved in the study of piece - wise linear homeomorphisms of the circle (i.e., orientation-preserving homeomorphisms with piece-wise constant derivatives) that if a piece-wise linear homeomorphism with an irrational rotation number \(\alpha\) is conjugate to an irrational rotation \({R}_{\alpha }\) through a \({C}^{1}\) diffeomorphism, then it is conjugate to an irrational rotation \({R}_{\alpha }\) through a piece-wise analytic homeomorphism. In this case, the unique invariant measure of \(f\) is equivalent to the Haar measure.

\section{Case of the Direct Product of Compact Groups}\

From the results of R. Arens, we know that for one-dimensional connected compact groups, the only measure-preserving homeomorphisms of the Haar measure are affine maps. In particular, for the unit circle, the only measure-preserving homeomorphisms are rotations and conjugate rotations. However, the generalization of R. Arens' results to higher dimensions is incorrect. The main work of this chapter is to find a general method for constructing measure-preserving homeomorphisms of the Haar measure on the direct product of compact groups. Through this method, we can construct non-affine measure-preserving homeomorphisms on the \(n\) -dimensional torus, which shows that Arens' results cannot be extended to higher-dimensional compact groups. At the same time, we apply this method to the torus and completely characterize the measure-preserving homeomorphisms of the Haar measure on the two-dimensional torus where one coordinate is a translation depending on the other coordinate. We uniformly use the following notations in this chapter and the next chapter. Let \(X\) be a compact set. Let \(\operatorname{Homeo}\left( X\right)\) denote the set of all homeomorphisms on \(X\). For any Borel measure \(\mu\) on \(X\), let \({\operatorname{Homeo}}_{\mu }\left( X\right)\) denote the set of all measure-preserving homeomorphisms of \(\mu\) on \(X\). For any compact groups \(G\) and \(x \in  G\) , let \({L}_{x}\) denote the left translation map on \(G\), that is, for any \(z \in  G,{L}_{x}\left( z\right)  = {xz}\). Let \(\operatorname{Aut}\left( G\right)\) denote the set of all maps on \(G\) that are both homeomorphisms and automorphisms. Let \(\operatorname{AF}\left( G\right)\) denote the set of all affine maps on \(G\), that is, \(\operatorname{AF}\left( G\right)  = \left\{  {{L}_{x} \circ  f : x \in  G,f \in  \operatorname{Aut}\left( G\right) }\right\}\).

\subsection{Construction of Measure-preserving Homeomorphisms}\
We introduce a method for constructing measure-preserving homeomorphisms of the Haar measure on \(n\) compact metric groups. In fact, we can construct measure-preserving homeomorphisms on \(n\) metric spaces with countable dense subsets. We have the following lemma. The main method of the proof is to use the integral characterization of measure-preserving maps of the Haar measure, that is, Theorem 2.1.5, and use Fubini's theorem to change the order of integration and handle each coordinate one by one.

\textbf{Lemma 2.1.1.} Let \({X}_{1},{X}_{2},\cdots ,{X}_{n}\) be a metric space with a countable dense subset, and \({\mu }_{1},{\mu }_{2},\cdots ,{\mu }_{n}\) be Borel measures on \({X}_{1},{X}_{2},\cdots ,{X}_{n}\) respectively. Let \({h}_{1}\) be a measure-preserving homeomorphism on \({X}_{1}\), and \({h}_{j}(2 \leq\) \(j \leq  n\)) be maps from \({X}_{1} \times  \cdots  \times  {X}_{j - 1}\) to \({\operatorname{Homeo}}_{{\mu }_{j}}\left( {X}_{j}\right)\) respectively. Let
$$
f ( x_{1}, \cdots, x_{n} )=\left( h_{1} ( x_{1} ), h_{2} \left( x_{1} \right) ( x_{2} ), \cdots, h_{n} ( x_{1}, \cdots, x_{n-1} ) \left( x_{n} \right) \right), 
$$
preserves the product measure \(\mu={\mu }_{1} \times  {\mu }_{2} \times  \cdots  \times  {\mu }_{n}\).

Proof: For any continuous function \(g\) on \({X}_{1} \times  {X}_{2} \times  \cdots  \times  {X}_{n}\), consider the integral of \(g\) on \({X}_{1} \times  {X}_{2} \times  \cdots  \times  {X}_{n}\)

By Fubini's theorem, we have
$$\begin{aligned}& \int_{X_1 \times X_2 \times \cdots \times X_n} g\left(f\left(x_1, x_2, \cdots, x_n\right)\right) d \mu \\
= & \int_{X_1 \times X_2 \times \cdots \times X_n} g\left(h_1\left(x_1\right), \cdots, h_n\left(x_1, \cdots, x_{n-1}\right)\left(x_n\right)\right) d \mu \\
= & \int_{X_1} \cdots \int_{X_n} g\left(h_1\left(x_1\right), \cdots, h_n\left(x_1, \cdots, x_{n-1}\right)\left(x_n\right)\right) d \mu \\
= & \int_{X_1} \cdots \int_{X_n} g\left(h_1\left(x_1\right), \cdots, h_{n-1}\left(x_1, \cdots, x_{n-2}\right)\left(x_{n-1}\right), x_n\right) d \mu \\
= & \int_{X_1} \cdots \int_{X_n} g\left(h_1\left(x_1\right), \cdots, h_{n-2}\left(x_1, \cdots, x_{n-3}\right)\left(x_{n-2}\right), x_{n-1}, x_n\right) d \mu \\
= & \cdots \cdots \cdots & \\
= & \int_{X_1} \cdots \int_{X_n} g\left(x_1, x_2, x_3, \cdots, x_{n-1}, x_n\right)d \mu.
\end{aligned}
$$
Therefore, \(f\left( {{x}_{1},{x}_{2},\cdots ,{x}_{n}}\right)\), \({\mu }_{1} \times  {\mu }_{2} \times  \cdots  \times  {\mu }_{n}\) preserves the product measure \({\mu }_{1} \times  {\mu }_{2} \times  \cdots  \times  {\mu }_{n}\).\hfill $\square$

From Lemma 2.1.1, it is easy to obtain the case of the direct product of \(n\) compact groups.

\textbf{Theorem 2.1.1.} Let \({G}_{1},{G}_{2},\cdots ,{G}_{n}\) be a compact metric group, and \({\mu }_{1},{\mu }_{2},\cdots ,{\mu }_{n}\) be the Haar measures on \({G}_{1},{G}_{2},\cdots ,{G}_{n}\) respectively. Let \(h\) be a self-homeomorphism on \({G}_{1}\) that preserves the Haar measure, and \({h}_{j}\left( {2 \leq  j \leq  n}\right)\) be the mappings from \({G}_{1} \times  \cdots  \times  {G}_{j - 1}\) to \({\operatorname{Homeo}}_{{\mu }_{j}}\left( {G}_{j}\right)\) respectively. Let
$$
f ( x_{1}, \cdots, x_{n} )=( h_{1} ( x_{1} ), h_{2} ( x_{1} ) ( x_{2} ), \cdots, h_{n} ( x_{1}, \cdots, x_{n-1} ) ( x_{n} ) ), 
$$
If \(f\left( {{x}_{1},{x}_{2},\cdots ,{x}_{n}}\right)\) is a self-homeomorphism on \({G}_{1} \times  {G}_{2} \times  \cdots  \times  {G}_{n}\), then \(f\left( {{x}_{1},{x}_{2},\cdots,{x}_{n}}\right)\) preserves the Haar measure \({\mu }_{1} \times  {\mu }_{2} \times  \cdots  \times  {\mu }_{n}\).

Using Theorem 2.1.1, we can easily construct some measure-preserving self-homeomorphisms of the product of compact groups. Below, we give some examples, and the most important one is the measure-preserving self-homeomorphism on the \(n\)-dimensional torus.

\textbf{Example 2.1.1.} Consider the \(n\)-dimensional torus \({\mathbb{T}}^{n} = \mathbb{T} \times  \cdots  \times  \mathbb{T} \cong  \mathbb{R}/\mathbb{Z} \times  \cdots  \times  \mathbb{R}/\mathbb{Z},\mathbb{T}\) as the unit circle in the complex plane. The topology on \(\mathbb{T}\) is the metric topology induced by the arc-length metric, and the topology on \({\mathbb{T}}^{n}\) is the product topology. The Lebesgue measure on \(\mathbb{T}\) is the Haar measure on \(\mathbb{T}\). The Haar measure on \({\mathbb{T}}^{n}\) is the product measure of the Haar measure on \(\mathbb{T}\). We use the functions and addition on \(\lbrack 0,1)\) to represent the functions and multiplication on the unit circle \(\mathbb{T}\). Let \(\alpha  \in  \lbrack 0,1),{g}_{j}\left( {2 \leq  j \leq  n}\right)\) be continuous functions from \({\mathbb{T}}^{j - 1} \cong  \mathbb{R}/\mathbb{Z} \times  \cdots  \times  \mathbb{R}/\mathbb{Z}\) to \(\mathbb{T} = \mathbb{R}/\mathbb{Z}\) respectively. Let \(f\left( {{x}_{1},{x}_{2},\cdots ,{x}_{n}}\right)  = \left( {{x}_{1} + \alpha ,{x}_{2} + {g}_{2}\left( {x}_{1}\right),\cdots ,{x}_{n} + {g}_{n}\left( {{x}_{1},\cdots ,{x}_{n - 1}}\right) }\right)\), then \(f\) is a Haar measure-preserving homeomorphism on \({G}_{1} \times  {G}_{2} \times  \cdots \times  {G}_{n}\).

Proof: Since \({g}_{j}\left( {2 \leq  j \leq  n}\right)\) is a continuous mapping, each coordinate component of \(f\) is continuous. Therefore, \(f\) is a continuous mapping. Next, we prove that \(f\) is a bijection. Suppose \(f\left( {{x}_{1},{x}_{2},\cdots ,{x}_{n}}\right)  = f\left( {{y}_{1},{y}_{2},\cdots ,{y}_{n}}\right)\). From the equality of the first coordinates, we can obtain \({x}_{1} = {y}_{1}\). From the equality of the second coordinates and \({x}_{1} = {x}_{2}\), we can obtain \({x}_{2} = {y}_{2}\). By recursively reasoning term by term, we can get \(\left( {{x}_{1},{x}_{2},\cdots ,{x}_{n}}\right)  = \left( {{y}_{1},{y}_{2},\cdots ,{y}_{n}}\right)\), that is, \(f\) is an injection. Let \(f\left( {{x}_{1},\cdots ,{x}_{n}}\right)  = \left( {{w}_{1},\cdots ,{w}_{n}}\right)\). From one coordinate, we can obtain \({x}_{1} = {w}_{1} - \alpha\). From the second coordinate and \({x}_{1} = {w}_{1} - \alpha\) , we can obtain \({x}_{2} = {w}_{2} - {g}_{2}\left( {{w}_{1} - }\right.\)  \(\alpha )\). By analogy, we can solve for \(\left( {{x}_{1},{x}_{2},\cdots ,{x}_{n}}\right)\). Therefore, \(f\) is a surjection. Then we have proved that \(f\) is a continuous bijection, so \(f\) is a self-homeomorphism on \({\mathbb{T}}^{n}\). In Theorem 2.1.1, take \({h}_{j}\left( {{x}_{1},\cdots ,{x}_{j - 1}}\right) \left( {x}_{j}\right)  =\)  \({x}_{j} + {g}_{j}\left( {{x}_{1},\cdots ,{x}_{j - 1}}\right) ,{h}_{1}\left( {x}_{1}\right)  = {x}_{1} + \alpha\). Since left-translation preserves the Haar measure, \(f\) preserves the Haar measure, that is, \(f\) is a measure-preserving self-homeomorphism on \({\mathbb{T}}^{n}\). \hfill $\square$

\textbf{Example 2.1.2.} Let the affine set of the unit circle \(\mathbb{T}\) be \(\mathrm{{AF}}\left( \mathbb{T}\right)\). We will prove in the next section that the affine transformation on a compact group preserves the Haar measure. Since the only maps on the circle that are both homeomorphisms and isomorphisms are the identity map and the inverse map, thus \(\mathrm{{AF}}\left( \mathbb{T}\right)  \cong  \mathbb{T} \rtimes  \{ 1, - 1\}\). Endow \(\{ 1, - 1\}\) with the discrete topology, then it is easy to prove that \(\mathrm{{AF}}\left( \mathbb{T}\right)\) is a compact group. We consider the product group of the compact groups \(\mathbb{T}\) and \(\mathrm{{AF}}\left( \mathbb{T}\right)\). Let \(g\) be a continuous function from \(\mathbb{T}\) to \(\mathbb{T}\), \(\alpha  \in  \lbrack 0,1),\beta  \in\)  \(\{ 1, - 1\}\), then \(f\left( {\left( {y,t}\right) ,x}\right)  = \left( {\left( {\alpha  + y,{\beta t}}\right) ,x + g\left( y\right) }\right)\) is a measure-preserving homeomorphism on \(\mathrm{{AF}}\left( \mathbb{T}\right)  \times  \mathbb{T}\).

Proof: In Theorem 2.1.1, we take \({h}_{1}\left( {y,t}\right)  = \left( {\alpha  + y,{\beta t}}\right) ,{h}_{2}\left( {y,t}\right)  = {L}_{g\left( y\right) }\). Since left-translation preserves the Haar measure, and since \(f\) is a homeomorphism on \(\mathrm{{AF}}\left( \mathbb{T}\right)  \times  \mathbb{T}\), so \(f\) preserves the Haar measure, that is, \(f\) is a measure-preserving homeomorphism on \(\mathrm{{AF}}\left( \mathbb{T}\right)  \times  \mathbb{T}\).\hfill $\square$

\subsection{Measure-preserving Homeomorphisms on the Torus}\

A special case of Example 2.1.1 in Section 2.1 on the two-dimensional torus \({\mathbb{T}}^{n}\) is an example of a self-homeomorphism that preserves the Haar measure and fixes one coordinate, that is, \(f\left( {x,y}\right)  = \left( {x,y + h\left( x\right) }\right)\), where \(h\) is a continuous function from \(\mathbb{T} = \mathbb{R}/\mathbb{Z}\) to \(\mathbb{T} = \mathbb{R}/\mathbb{Z}\). The main work in this section is to give a complete characterization of such measure-preserving self-homeomorphisms that fix one coordinate and measure-preserving self-homeomorphisms where the first coordinate is a translation depending on the second coordinate, and generalize it to the \(n\)-dimensional torus. We need two lemmas before giving the characterization.

The following lemma comes from the result of R. Arens, that is, all self-homeomorphisms of one-dimensional connected compact topological groups that preserve the Haar measure are affine. In particular, all self-homeomorphisms of the circle that preserve the Haar measure are affine. However, for this specific case, we can prove it by a simple method.

\textbf{Lemma 2.2.1.} Let \(f\) be a homeomorphism on the unit circle \(\mathbb{T} = \mathbb{R}/\mathbb{Z}\) that preserves the Haar measure. Then there exists \(\alpha  \in  \lbrack 0,1)\) such that \(f\left( x\right)  = \alpha  + x\) or \(f\left( x\right)  = \alpha  - x\).

Proof: Let \(g\left( x\right)  = f\left( x\right)  - f\left( 0\right)\). Then \(g\left( 0\right)  = 0\left( {\;\operatorname{mod}\;1}\right)\). Since translation preserves the Haar measure, \(g\) preserves the Haar measure. Since \(g\) is a homeomorphism on \(\lbrack 0,1)\), \(g\) is strictly monotonic on \(\lbrack 0,1)\). Also, since \(g\left( 0\right)  = 0\left( {\;\operatorname{mod}\;1}\right)\), then \(g\left( 0\right)  = 0\) or \(g\left( 0\right)  = 1\). If \(g\left( 0\right)  = 0\), then \(g\) is strictly increasing on \(\lbrack 0,1)\). Therefore, \(g\left( 1\right)  = 1\). Since \(g\) preserves the Haar measure \(\mu\), for any \(x \in  \lbrack 0,1),\mu \left( {g\left( \left\lbrack  {0,x}\right\rbrack  \right) }\right)  =\)  \(\mu \left( \left\lbrack  {0,x}\right\rbrack  \right)\), that is, \(\mu \left( \left\lbrack  {g\left( 0\right) ,g\left( x\right) }\right\rbrack  \right)  = \mu \left( \left\lbrack  {0,x}\right\rbrack  \right)\), then \(g\left( x\right)  - g\left( 0\right)  = x\). Therefore, \(g\left( x\right)  = x\). So \(f\left( x\right)  =\)  \(f\left( 0\right)  + x\). If \(g\left( 0\right)  = 1\), then \(g\) is strictly decreasing on \(\lbrack 0,1)\). Therefore, \(g\left( 1\right)  = 1\). Since \(g\) preserves the Haar measure \(\mu\), for any \(x \in  \lbrack 0,1),\mu \left( {g\left( \left\lbrack  {0,x}\right\rbrack  \right) }\right)  = \mu \left( \left\lbrack  {0,x}\right\rbrack  \right)\), that is, \(\mu \left( \left\lbrack  {g\left( 0\right) ,g\left( x\right) }\right\rbrack  \right)  = \mu \left( \left\lbrack  {0,x}\right\rbrack  \right)\), then \(g\left( x\right)  =  - x\). Therefore, \(f\left( x\right)  = f\left( 0\right)  - x\). \hfill $\square$

\textbf{Lemma 2.2.2.} Let \(f\left( {x,y}\right)  = \left( {x,g\left( {x,y}\right) }\right)\) be a homeomorphism on the two-dimensional torus \({\mathbb{T}}^{2}\) that preserves the Haar measure. Then for any fixed \({x}_{0} \in  \lbrack {0.1})\left( {\;\operatorname{mod}\;1}\right) ,g\left( {{x}_{0},y}\right)\) is a homeomorphism on the unit circle \(y\) that preserves the Haar measure.

Proof: Denote the Haar measure on the unit circle \(\mathbb{T}\) as \(m\), and the Haar measure on the two-dimensional torus \({\mathbb{T}}^{2}\) as \(\mu\). For any Borel set \(B\) on the unit circle \(\mathbb{T}\), for any \({x}_{0} \in  \lbrack {0.1})\left( {\;\operatorname{mod}\;1}\right)\) and \(\varepsilon  > 0\), we have
$$
\left(x_0-\varepsilon, x_0+\varepsilon\right) \times\left\{g\left(x_0, y\right): y \in B\right\} \subseteq f\left(\left(x_0-\varepsilon, x_0+\varepsilon\right) \times B\right)
$$
and
$$
f\left(\left(x_0-\varepsilon, x_0+\varepsilon\right) \times B\right) \subseteq E
$$
Where $E=\left( {{x}_{0} - \varepsilon ,{x}_{0} + \varepsilon }\right)  \times  \left\{  {g\left( {x,y}\right)  : x \in  \left( {{x}_{0} - \varepsilon ,{x}_{0} + \varepsilon }\right) ,y \in  B}\right\}$.
Denote $C = \left\{  {g\left( {{x}_{0},y}\right)  : y \in  B}\right\}$, and
$${D}_{\varepsilon } = \left\{  {g\left( {x,y}\right)  : x \in  \left( {{x}_{0} - \varepsilon ,{x}_{0} + \varepsilon }\right) ,y \in  B}\right\},$$ then  
$$
\mu \left( {\left( {{x}_{0} - \varepsilon ,{x}_{0} + \varepsilon }\right)  \times  C}\right)  \leq  \mu \left( {f\left( {\left( {{x}_{0} - \varepsilon ,{x}_{0} + \varepsilon }\right)  \times  B}\right) }\right).
$$
$$
\mu \left( {f\left( {\left( {{x}_{0} - \varepsilon ,{x}_{0} + \varepsilon }\right)  \times  B}\right) }\right)  \leq  \mu \left( {\left( {{x}_{0} - \varepsilon ,{x}_{0} + \varepsilon }\right)  \times  {D}_{\varepsilon }}\right).
$$
Since \(f\) is a homeomorphism on the two-dimensional torus \({\mathbb{T}}^{2}\) that preserves the Haar measure, then
\[
\mu \left( {f\left( {\left( {{x}_{0} - \varepsilon ,{x}_{0} + \varepsilon }\right)  \times  B}\right) }\right)  = \mu \left( {\left( {{x}_{0} - \varepsilon ,{x}_{0} + \varepsilon }\right)  \times  B}\right),
\]
Therefore
${2\varepsilon } \times  m\left( C\right)  \leq  {2\varepsilon } \times  m\left( B\right)  \leq  {2\varepsilon } \times  m\left( {D}_{\varepsilon }\right)$,
Therefore
\[
m\left( C\right)  \leq  m\left( B\right)  \leq  m\left( {D}_{\varepsilon }\right).
\]
Take \(\varepsilon  = \frac{1}{n}\), denote \({D}_{n} = \left\{  {g\left( {x,y}\right)  : x \in  \left( {{x}_{0} - \frac{1}{n},{x}_{0} + \frac{1}{n}}\right) ,y \in  B}\right\}\), then \({D}_{1} \supseteq  {D}_{2} \supseteq  \cdots  \supseteq  C\), and \(\mathop{\bigcap }\limits_{{n = 1}}^{\infty }{D}_{n} = C\).
Therefore \(\mathop{\lim }\limits_{{n \rightarrow  \infty }}m\left( {D}_{n}\right)  = m\left( C\right)\), and \(m\left( B\right)\) and \(m\left( C\right)\) are independent of \(n\). 
Therefore \(m\left( C\right)  = m\left( B\right)\), 
that is, 
$
2 \varepsilon \times m(C) \leq 2 \varepsilon \times m(B) \leq 2 \varepsilon \times m\left(D_{\varepsilon}\right),
$
For any ${x}_{0} \in \mathbb{T}$, $g\left( {{x}_{0},y}\right)$ is a continuous bijection on the unit circle \(\mathbb{T}\) . Since \(\mathbb{T}\) is a compact group, \(g\left( {{x}_{0},y}\right)\) is a homeomorphism. This shows that \(g\left( {{x}_{0},y}\right)\) is a homeomorphism on the unit circle \(\mathbb{T}\) that preserves the Haar measure with respect to \(y\). \hfill $\square$

\textbf{Theorem 2.2.1.} Let \(f\left( {x,y}\right)  = \left( {x,g\left( {x,y}\right) }\right)\) be a homeomorphism on the two-dimensional torus \({\mathbb{T}}^{2}\) that preserves the Haar measure. Then there exists a continuous function \(h\) from \(\mathbb{T}\) to \(\mathbb{T}\) such that \(f\left( {x,y}\right)  = \left( {x,h\left( x\right)  + y}\right)\) or \(f\left( {x,y}\right)  = \left( {x,h\left( x\right)  - y}\right)\).

Proof: Since \(f\left( {x,y}\right)  = \left( {x,g\left( {x,y}\right) }\right)\) is a homeomorphism on the two-dimensional torus \({\mathbb{T}}^{2}\) that preserves the Haar measure, by Lemma 2.2.2, for any \(x \in  \lbrack {0.1})\left( {\;\operatorname{mod}\;1}\right),g\left( {x,y}\right)\) is a measure-preserving homeomorphism with respect to \(y\). By Lemma 2.2.1, there exist a number \(\alpha \left( x\right)  \in  \lbrack 0,1)\) related to \(x\) and \(\beta \left( x\right)  \in  \{ 1, - 1\}\) such that \(g\left( {x,y}\right)  = \alpha \left( x\right)  + \beta \left( x\right) y\). Since \(g\left( {x,y}\right)\) is jointly continuous with respect to \(x,y\), then for a fixed \(y \in  \lbrack {0.1})\left( {\;\operatorname{mod}\;1}\right), g\left( {x,y}\right)  = \alpha \left( x\right)  + \beta \left( x\right) y\), it is continuous with respect to \(x\). Taking \(y = 0\), then \(g\left( {x,0}\right)  = \alpha \left( x\right)\) is continuous with respect to \(y\). Therefore, \(\beta \left( x\right) y\) is jointly continuous with respect to \(x,y\). Then for any \(y \in  \mathbb{T}\), that is, \({e}^{{2\pi i\beta }\left( x\right) y}\) is continuous with respect to \(x\), and also for any \(z \in  \mathbb{T},{z}^{\beta \left( x\right) }\) is continuous with respect to \(x\). Let \(t = {e}^{2\pi ix}\), then \(\beta \left( x\right)  = \gamma \left( t\right)\) is a continuous function of \(t \in  \mathbb{T}\). For any \(z \in  \mathbb{T},{z}^{\gamma \left( t\right) }\), it is continuous with respect to \(t\). Next, we prove that for any \(t \in  \mathbb{T},\gamma \left( t\right)\) is a constant. We use the method of contradiction. Assume that \(t \in  \mathbb{T}\) is not a constant. First, we prove that there exists \({t}_{0} \in  \mathbb{T}\) such that for any neighborhood \(U,{z}^{\gamma \left( t\right) }\) of \({t}_{0}\), it is not a constant on \(U\). Otherwise, for any \(t \in  \mathbb{T}\), there exists a neighborhood \({U}_{t}\) of \(t\) such that \(\gamma \left( t\right)\) is a constant on \({U}_{t}\). For each \({U}_{t}\), take \(\delta\) small enough such that \({O}_{t} = \{ s \in  \mathbb{T} : d\left( {s,t}\right)  < \delta \}  \subseteq  {U}_{t}\), where \(d\) is the arc length metric. Then \(\gamma \left( t\right)\) is a constant on \({O}_{t}\). Since \(\mathop{\bigcup }\limits_{{t \in  {S}^{1}}}{O}_{t} = \mathbb{T}\) and the unit circle \(\mathbb{T}\) is a compact group, there exists \({t}_{1},{t}_{2},\cdots ,{t}_{n} \in  \mathbb{T}\) such that \(\mathop{\bigcup }\limits_{{i = 1}}^{n}{O}_{{t}_{i}} = \mathbb{T}\). Since \(\mathbb{T}\) is a connected set, for \({O}_{{i}_{1}}\), there must exist \({O}_{{i}_{2}}\) such that \({O}_{{i}_{1}} \cap  {O}_{{i}_{2}}\) is non-empty. Therefore, \(\gamma \left( t\right)\) is a constant on \({O}_{{i}_{1}} \cup  {O}_{{i}_{2}}\). For \({O}_{{i}_{1}} \cup  {O}_{{i}_{2}}\), there must exist \({O}_{{i}_{3}}\) such that \(\left( {{O}_{{i}_{1}} \cup  {O}_{{i}_{2}}}\right)  \cap  {O}_{{i}_{3}}\) is non-empty. Therefore, \(\gamma \left( t\right)\) is a constant on \({O}_{{i}_{1}} \cup  {O}_{{i}_{2}} \cup  {O}_{{i}_{3}}\). Repeating this process, since \({t}_{1},{t}_{2},\cdots ,{t}_{n} \in  \mathbb{T}\) is finite, it will stop after a certain step. Then \(\gamma \left( t\right)\) is a constant on \(\mathbb{T}\), which contradicts the fact that \(\gamma \left( t\right)\) is not a constant function on \(\mathbb{T}\). Therefore, there exists \({t}_{0} \in  \mathbb{T}\) such that for any neighborhood \(U,\gamma \left( t\right)\) of \({t}_{0}\), it is not a constant on \(U\). Also, since for any \(z \in  \mathbb{T},{z}^{\gamma \left( t\right) }\) is continuous with respect to \(t\), when \(t \rightarrow  {t}_{0}\), \({z}^{\gamma \left( t\right) } \rightarrow  {z}^{\gamma \left( {t}_{0}\right) }\). Therefore, there exists an integer sequence \({\left\{  {k}_{n}\right\}  }_{n = 1}^{\infty }\), where \({k}_{n} \in  \{  - 2, - 1,1,2\}\), such that for any \(z \in  \mathbb{T}\), \({z}^{{k}_{n}} \rightarrow  1\left( {n \rightarrow  \infty }\right)\). Next, we prove \(\mu \left( \left\{  {z \in  \mathbb{T} : {z}^{{k}_{n}} \rightarrow  1,n \rightarrow  \infty }\right\}  \right)  = 0\) to get a contradiction. Consider the set
$$
\left\{z \in \mathbb{T}: z^{k_n} \rightarrow 1, n \rightarrow \infty\right\}=\bigcap_{s=2}^{\infty} \bigcup_{N=1}^{\infty} \bigcap_{n=N}^{\infty}\left\{z \in \mathbb{T}: d\left(z^{k_n}-1\right) \leq \frac{1}{s}\right\},
$$
For any \(n \geq  N,\mu \left( \left\{  {z \in  \mathbb{T} : d\left( {{z}^{{k}_{n}} - 1}\right)  \leq  \frac{2}{s}}\right\}  \right)  = \frac{2}{s}\) , therefore
$$
\begin{aligned}
\mu\left(\left\{z \in \mathbb{T}: z^{k_n} \rightarrow 1\right\}\right) & =\mu\left(\bigcap_{s=2}^{\infty} \bigcup_{N=1}^{\infty} \bigcap_{n=N}^{\infty}\left\{z \in \mathbb{T}: d\left(z^{k_n}-1\right) \leq \frac{1}{s}\right\}\right) \\
& \leq \mu\left(\bigcap_{s=2}^{\infty} \bigcup_{N=1}^{\infty}\left\{z \in \mathbb{T}: d\left(z^{k_N}-1\right) \leq \frac{1}{s}\right\}\right) \\
& =\lim _{s \rightarrow \infty} \frac{2}{s}=0.
\end{aligned}
$$
However, for any \(z \in  \mathbb{T}\), there is always \({z}^{{k}_{n}} \rightarrow  1\left( {n \rightarrow  \infty }\right)\). Therefore, \(\left\{  {z \in  \mathbb{T} : {z}^{{k}_{n}} \rightarrow  1,n \rightarrow  \infty }\right\}   = \mathbb{T}\), and \(\mu \left( \left\{  {z \in  \mathbb{T} : {z}^{{k}_{n}} \rightarrow  1,n \rightarrow  \infty }\right\}  \right)  = 1\), a contradiction! This shows that for any \(t \in  \mathbb{T},\gamma \left( t\right)\), it is a constant, that is, \(\beta \left( x\right)\) is a constant. Therefore, \(g\left( {x,y}\right)  = \alpha \left( x\right)  + y\) or \(g\left( {x,y}\right)  = \alpha \left( x\right)  - y\). Among them, \(\alpha \left( x\right)\) is a continuous function from \(\mathbb{T}\) to \(\mathbb{T}\). \hfill $\square$

Since \(g\left( {x,y}\right)  = \left( {y,x}\right) ,g\left( {x,y}\right)  = \left( {x, - y}\right)\) are all affine transformations, Theorem 2.2.1 shows that a measure-preserving self-homeomorphism with a fixed coordinate must be a composition of a measure-preserving homeomorphism of the form \(f\left( {x,y}\right)  = \left( {x,h\left( x\right)  + y}\right)\) and an affine transformation. We can generalize this result to measure-preserving homeomorphisms of higher-dimensional tori.

\textbf{Proposition 2.2.1.} Let
$$
f\left( {{x}_{1},\cdots ,{x}_{n}}\right)  = \left( {{h}_{1}\left( {x}_{1}\right) ,{h}_{2}\left( {{x}_{1},{x}_{2}}\right) ,\cdots ,{h}_{n}\left( {{x}_{1},{x}_{2},\cdots ,{x}_{n}}\right) }\right) 
$$
 be a self-homeomorphism on the \(n\)-dimensional torus \({\mathbb{T}}^{n}\) that preserves the Haar measure, where \({h}_{1}\) is a measure-preserving self-homeomorphism on the unit circle \(\mathbb{T}\) . Then there exist \(\alpha  \in  \lbrack 0,1)\) and continuous functions \({g}_{1},\cdots ,{g}_{n - 1}\) from \(\mathbb{T},{\mathbb{T}}^{2},\cdots ,{\mathbb{T}}^{n - 1}\) to \(\mathbb{T}\) such that \(f\left( {{x}_{1},\cdots ,{x}_{n}}\right)  = \left( {\pm {x}_{1} + \alpha , \pm  {x}_{2} + {g}_{1}\left( {x}_{1}\right) ,\cdots , \pm  {x}_{n} + {g}_{n - 1}\left( {{x}_{1},\cdots ,{x}_{n - 1}}\right) }\right)\).

Proof: By induction, for the case of \(n = 2\), it follows from Theorem 2.2.1 and Lemma 2.2.2. Assume that the case of \(n - 1\) dimensions holds. For the case of \(n\) dimensions, using the same proof method as in Lemma 2.2.2, it is easy to prove that for any fixed \({x}_{1},{x}_{2},\cdots ,{x}_{n - 1} \in  \mathbb{T},{h}_{n}\left( {{x}_{1},{x}_{2},\cdots ,{x}_{n}}\right)\), it is a measure-preserving self-homeomorphism with respect to \({x}_{n}\). Similar to the proof of Theorem 2.2.1, it can be proved that there exists a continuous function \({g}_{n - 1}\) from \({\mathbb{T}}^{n - 1}\) to \(\mathbb{T}\) such that \({h}_{n}\left( {{x}_{1},{x}_{2},\cdots ,{x}_{n}}\right)  =  \pm  {x}_{n} + {g}_{n - 1}\left( {{x}_{1},\cdots ,{x}_{n - 1}}\right)\). Let
\[
k\left( {{x}_{1},{x}_{2},\cdots ,{x}_{n}}\right)  = \left( {{x}_{1},{x}_{2},\cdots ,{x}_{n - 1}, \pm  {x}_{n} + {g}_{n - 1}\left( {{x}_{1},\cdots ,{x}_{n - 1}}\right) }\right)
\]
and
\[
l\left( {{x}_{1},{x}_{2},\cdots ,{x}_{n}}\right)  = \left( {{h}_{1}\left( {x}_{1}\right) ,{h}_{2}\left( {{x}_{1},{x}_{2}}\right) ,\cdots ,{h}_{n - 1}\left( {{x}_{1},{x}_{2},\cdots ,{x}_{n - 1}}\right) ,{x}_{n}}\right),
\]
\[
{l}_{1}\left( {{x}_{1},{x}_{2},\cdots ,{x}_{n - 1}}\right)  = \left( {{h}_{1}\left( {x}_{1}\right) ,{h}_{2}\left( {{x}_{1},{x}_{2}}\right) ,\cdots ,{h}_{n - 1}\left( {{x}_{1},{x}_{2},\cdots ,{x}_{n - 1}}\right) }\right).
\]
Then \(f = l \circ  k\) . Since the first \(n - 1\) coordinates of \(l\left( {{x}_{1},\cdots ,{x}_{n}}\right)  = \left( {{h}_{1}\left( {x}_{1}\right) ,{h}_{2}\left( {{x}_{1},{x}_{2}}\right) ,\cdots ,{h}_{n - 1}\left( {{x}_{1},{x}_{2},\cdots ,{x}_{n - 1}}\right) ,{x}_{n}}\right)\) are independent of \({x}_{n}\) , \({l}_{1}\left( {{x}_{1},{x}_{2},\cdots ,{x}_{n - 1}}\right)\) is a self-homeomorphism on the \(n - 1\)-dimensional torus \({\mathbb{T}}^{n - 1}\) that preserves the Haar measure. By the inductive hypothesis, there exist \(\alpha  \in  \lbrack 0,1)\) and continuous functions \({g}_{1},\cdots ,{g}_{n - 2}\) from \(\mathbb{T},{\mathbb{T}}^{2},\cdots ,{\mathbb{T}}^{n - 2}\) to \(\mathbb{T}\) respectively such that
$${l}_{1}\left( {{x}_{1},\cdots ,{x}_{n - 1}}\right)  = \left( {\pm {x}_{1} + \alpha ,\cdots , \pm  {x}_{n - 1} + {g}_{n - 2}\left( {{x}_{1},\cdots ,{x}_{n - 2}}\right) }\right),$$
Therefore 
$$
f\left( {{x}_{1},\cdots ,{x}_{n}}\right)  = \left( {\pm {x}_{1} + \alpha , \pm  {x}_{2} + {g}_{1}\left( {x}_{1}\right) ,\cdots , \pm  {x}_{n} + {g}_{n - 1}\left( {{x}_{1},\cdots ,{x}_{n - 1}}\right) }\right).
$$ \hfill $\square$

Furthermore, for the measure-preserving self-homeomorphism on the two-dimensional torus where the first coordinate depends on the translation of the second coordinate, we can also obtain its complete characterization.

\textbf{Theorem 2.2.2.} Let \(f\left( {x,y}\right)  = \left( {x + {g}_{1}\left( y\right) ,h\left( {x,y}\right) }\right)\) be a homeomorphism on the two-dimensional torus \({\mathbb{T}}^{2}\) that preserves the Haar measure, where \({g}_{1}\) is a continuous function from \(\mathbb{T}\) to \(\mathbb{T}\) . Then there exists a continuous function \({g}_{2}\) from \(\mathbb{T}\) to \(\mathbb{T}\) such that 
$$f\left( {x,y}\right)  = \left( {x + {g}_{1}\left( y\right) ,{g}_{2}\left( {x + {g}_{1}\left( y\right) }\right)  + y}\right)$$
or 
$$f\left( {x,y}\right)  = \left( {x + {g}_{1}\left( y\right) ,{g}_{2}\left( {x + {g}_{1}\left( y\right) }\right)  - y}\right).$$

Proof: For any \(\left( {x,y}\right)  \in  {\mathbb{T}}^{2}\) , let \(k\left( {x,y}\right)  = \left( {x - {g}_{1}\left( y\right) ,y}\right)\). Then, according to Example \({2.1.1},k\), \({\mathbb{T}}^{2}\) is a measure-preserving homeomorphism, and \(f \circ  k\left( {x,y}\right)  = \left( {x,h\left( {x - {g}_{1}\left( y\right) ,y}\right) }\right)\). By Theorem 2.2.1, there exists a continuous function \({g}_{2}\) from \(\mathbb{T}\) to \(\mathbb{T}\) such that \(f \circ  k\left( {x,y}\right)  = \left( {x + {g}_{1}\left( y\right) ,{g}_{2}\left( x\right)  + y}\right)\) or \(f \circ  k\left( {x,y}\right)  = \left( {x + {g}_{1}\left( y\right) ,{g}_{2}\left( x\right)  - y}\right)\). Then 
$$f\left( {x,y}\right)  = \left( {x + {g}_{1}\left( y\right) ,{g}_{2}\left( {x + {g}_{1}\left( y\right) }\right)  + y}\right)$$
or
$$f\left( {x,y}\right)  = \left( {x + {g}_{1}\left( y\right) ,{g}_{2}\left( {x + {g}_{1}\left( y\right) }\right)  - y}\right).$$ \hfill $\square$

A special case of this theorem is \(f\left( {x,y}\right)  = \left( {x + y,h\left( {x,y}\right) }\right)\). In this case, there exists a continuous function \(g\) from \(\mathbb{T}\) to \(\mathbb{T}\) such that \(f\left( {x,y}\right)  = \left( {x + y,g\left( {x + y}\right)  + y}\right)\) or \(f\left( {x,y}\right)  = \left( {x + y,g\left( {x + y}\right)  - y}\right)\). We can use these results to determine whether some specific homeomorphisms are measure-preserving homeomorphisms. For example, let \(f\left( {x,y}\right)  =\)  \(\left( {x + y,x + \sin \left( {2\pi y}\right) /{2\pi }}\right)\). It can be proved that it is a homeomorphism on the two-dimensional torus, but it does not have the above form. Therefore, it is easy to prove by contradiction that it is not a homeomorphism that preserves the Haar measure.

\section{Generalization and Application of the Normalizer}\

In this chapter, we defined the concept of the generalized normalizer of a subgroup and proved that under certain conditions, the generalized normalizer is a set composed of measure-preserving homeomorphisms. In particular, the mappings in the normalizer of an affine set must preserve the Haar measure. In the first section, we showed that in some non-commutative compact groups, the normalizer of an affine set is strictly larger than the affine set, which indicates that new measure-preserving homeomorphisms can indeed be obtained through the normalizer. And since the proof is constructive, we also obtained a new method for constructing measure-preserving homeomorphisms. However, not all normalizer subsets of affine sets in compact groups can expand the affine set. In the second section, we gave an example of a compact group where the affine set is strictly equal to its normalizer.

\subsection{Properties and Applications of the Normalizer Subset}\

In this section, we defined the concept of the generalized normalizer on a compact group. Through this concept, we can obtain new measure-preserving homeomorphisms from the known set of measure-preserving homeomorphisms, thus getting a new method for constructing measure-preserving homeomorphisms. Not only for the Haar measure in compact groups, we can define the concept of the normalizer for general Borel measures on topological spaces.

\textbf{Definition 3.1.1.} Let \(X\) be a topological space, \(\mu\) be a Borel probability measure on \(X\), \(K,P \subseteq  {\operatorname{Homeo}}_{\mu }\left( X\right)\). Define the normalizer from \(K\) to \(P\) as
\[
{E}_{K}\left( P\right)  = \{ f \in  \operatorname{Homeo}\left( X\right)  : {fK} \subseteq  {Pf}\}.
\]
\textbf{Proposition 3.1.1.} Let \(P,{P}_{1},{P}_{2},K,{K}_{1},{K}_{2} \subseteq  {\operatorname{Homeo}}_{\mu }\left( X\right) ,{K}_{2} \subseteq  {K}_{1},{P}_{1} \subseteq  {P}_{2}\), then

(1) If \({K}_{1} \subseteq  {K}_{2}\), then \({E}_{{K}_{1}}\left( P\right)  \supseteq  {E}_{{K}_{2}}\left( P\right)\);

(2) If \({P}_{1} \subseteq  {P}_{2}\), then \({E}_{K}\left( {P}_{1}\right)  \subseteq  {E}_{K}\left( {P}_{2}\right)\).

Proof: If \({K}_{1} \subseteq  {K}_{2}\), for any \(f \in  {E}_{{K}_{2}}\left( P\right)\) , we have \(f{K}_{2} \subseteq  {Pf}\), then
$
f{K}_{1} \subseteq  f{K}_{2} \subseteq  {Pf}.
$
Therefore \(f \in  {E}_{{K}_{1}}\left( P\right)\), that is \({E}_{{K}_{1}}\left( P\right)  \supseteq  {E}_{{K}_{2}}\left( P\right)\).

If \({P}_{1} \subseteq  {P}_{2}\), for any \(f \in  {E}_{K}\left( {P}_{1}\right)\), we have \({fK} \subseteq  {P}_{1}f\), then
$
{fK} \subseteq  {P}_{1}f \subseteq  {P}_{2}f.
$
Therefore \({E}_{K}\left( {P}_{1}\right)  \subseteq  {E}_{K}\left( {P}_{2}\right)\). \hfill $\square$

Proposition 3.1.1 shows the monotonicity of the normalizer from \(K\) to \(P\) with respect to \(K\) and \(P\). We will use this property several times. Next, we prove that if \(K\) has a unique invariant probability measure, then the elements in the normalizer are also measure-preserving homeomorphisms.

\textbf{Theorem 3.1.1.} Let \(X\) be a compact topological space, \(\mu\) be a Borel probability measure on \(X\), \(P,K \subseteq  {\operatorname{Homeo}}_{\mu }\left( X\right)\) and \(\mu\) be the unique measure that is an invariant Borel probability measure for all elements in \(K\), then \({E}_{K}\left( P\right)  \subseteq\)  \({\operatorname{Homeo}}_{\mu }\left( X\right)\).

Proof: For any Borel set \(B\), any \(f \in  {E}_{K}\left( P\right)\), let \({\mu }_{f}\left( B\right)  = \mu \left( {f\left( B\right) }\right)\). Next, we prove that \({\mu }_{f}\) is
a Borel probability measure. For any sequence of pairwise-disjoint Borel sets \({\left\{  {A}_{i}\right\}  }_{i = 1}^{\infty }\), since \(f\) is a homeomorphism, then
$$
\begin{aligned}
\mu_f\left(\bigcup_{i=1}^{\infty} A_i\right) & =\mu\left(f\left(\bigcup_{i=1}^{\infty} A_i\right)\right) \\
& =\mu\left(\bigcup_{i=1}^{\infty} f\left(A_i\right)\right) \\
& =\sum_{i=1}^{\infty} \mu_f\left(A_i\right).
\end{aligned}
$$
Therefore \({\mu }_{f}\) satisfies countable additivity, and since
$$
\mu_f(X)=\mu(f(X))=\mu(f(X))=1
$$
Therefore, \({\mu }_{f}\) is indeed a Borel probability measure. Since \(f \in  {E}_{K}\left( P\right)\), for any \(k \in  K\), there exists \(p \in  P\) such that \(f \circ  k = p \circ  f\). Therefore, for any Borel set \(B\), since \(p \in  P\) preserves the measure \(\mu\), thus
\[
{\mu }_{f}(k\left( B\right)  = \mu \left( {f \circ  k\left( B\right) }\right)  = \mu \left( {p \circ  f\left( B\right) }\right)  = \mu \left( {f\left( B\right) }\right)  = {\mu }_{f}\left( B\right),
\]
Since \(\mu\) is the unique \(K\)-invariant measure, then \({\mu }_{f} = \mu\), that is, for any Borel set \(B,\mu \left( {f\left( \mathrm{\;B}\right) }\right)  =\)  \(\mu \left( B\right)\). Because \(f\) is a homeomorphism, \(B\) is a Borel set if and only if there exists a Borel set \(A\) such that \(f\left( A\right)  = B\). Therefore
\[
\mu \left( {{f}^{-1}\left( B\right) }\right)  = \mu \left( A\right)  = \mu \left( {f\left( A\right) }\right)  = \mu \left( B\right).
\]
So \(f \in  {E}_{K}\left( P\right)\) preserves the measure, that is, \({E}_{K}\left( P\right)  \subseteq  {\operatorname{Homeo}}_{\mu }\left( X\right)\). \hfill $\square$

Theorem 3.1.1 still holds when \(\mu\) is a regular Borel probability measure and \(\mu\) is the unique invariant regular Borel probability measure for \(K\) .

\textbf{Corollary 3.1.1.} Let \(X\) be a compact topological space, \(\mu\) be a regular Borel probability measure on \(X\) , \(P,K \subseteq\)  \({\text{ Homeo }}_{\mu }\left( X\right)\) and \(\mu\) be the unique measure that is a regular invariant Borel probability measure for all elements in \(K\) . Then \({E}_{K}\left( P\right)  \subseteq  {\operatorname{Homeo}}_{\mu }\left( X\right)\) .

Proof: For any Borel set \(B\) and any \(f \in  E\left( P\right)\) , let \({\mu }_{f}\left( B\right)  = \mu \left( {f\left( B\right) }\right)\) . From the proof of Lemma 3.2.1, it suffices to prove that \({\mu }_{f}\) is a regular measure. Since \(f\) is a self-homeomorphism of \(G\) , any open set \(U \subseteq  G\) if and only if there exists an open set \(O \subseteq  G\) , and any compact set \(F \subseteq  G\) if and only if there exists a compact set \(J \subseteq  G\) . Therefore
$$
\begin{aligned}
\mu_f(B) & =\mu(f(B)) \\
& =\inf \{\mu(U): f(B) \subseteq U, U \in G \text { is open set }\} \\
& =\inf \{\mu(f(O)): f(B) \subseteq f(O), O \in G \text { is open set }\} \\
& =\inf \left\{\mu_f(O): \mathrm{B} \subseteq O, O \in G \text { is open set }\right\},\\
\mu_f(B) & =\mu(f(B)) \\
& =\sup \{\mu(F): F \subseteq f(B), F \in G \text { is compact set }\} \\
& =\sup \{\mu(f(J)): f(J) \subseteq f(B), J \in G \text { is compact set }\} \\
& =\sup \left\{\mu_f(J): J \subseteq B, J \in G \text { is compact set }\right\}.
\end{aligned}
$$
So \({\mu }_{f}\) also has regularity. Therefore, \(f\) preserves the Haar measure \(\mu\), that is, \({E}_{K}\left( P\right)  \subseteq  {\operatorname{Homeo}}_{\mu }\left( G\right)\). \hfill $\square$

Theorem 3.1.1 obtains a new set of measure-preserving homeomorphisms from a known set of measure-preserving self-homeomorphisms. Therefore, as long as we can find elements in \({E}_{K}\left( P\right)  \smallsetminus  \left( {P \cup  K}\right)\), we get new measure-preserving self-homeomorphisms.

Since the Haar measure plays a crucial role in harmonic analysis and serves as the foundation of the entire modern harmonic analysis, and it is closely related to the maximum entropy of the dynamical system and also plays an important role in the dynamical system, being the most important measure on a compact group, we mainly consider the case where \(X\) is a compact group \(G\) and \(\mu\) is the Haar measure. Since the Haar measure is the unique left-translation invariant regular Borel probability measure on a compact group, the left-translation set \(\left\{  {{L}_{x} : x \in  G}\right\}\) is an important set of measure-preserving self-homeomorphisms on \(G\) . Therefore, in many cases, we take \(K =\)  \(\left\{  {{L}_{x} : x \in  G}\right\}\). Unless otherwise emphasized, we default that the normalizer of the set \(P\) of measure-preserving self-homeomorphisms on the compact group \(G\) is from \(K = \left\{  {{L}_{x} : x \in  G}\right\}\) to \(P\). We abbreviate \({E}_{\left\{  {L}_{x} : x \in  G\right\}  }\left( P\right)\) as \({E}_{G}\left( P\right)\), which is called the translation normalizer of \(P\) on \(G\) . By Corollary 3.2.1, \({E}_{G}\left( P\right)  \subseteq  {\operatorname{Homeo}}_{\mu }\left( G\right)\).

Since \(G\) has a group structure, we are more concerned about when \({E}_{K}\left( P\right)\) has a group structure. The following proposition gives a sufficient condition for \({E}_{K}\left( P\right)\) to be a group.

\textbf{Proposition 3.1.2.} Let \(G\) be a compact group and \(\mu\) be the Haar measure on \(G\). If \(P \subseteq  K \subseteq  {\operatorname{Homeo}}_{\mu }\left( G\right)\), then \({E}_{K}\left( P\right)\) is a group under the composition operation of mappings.

Proof: For any \(f,g \in  {E}_{K}\left( P\right)\), then for any \(k \in  K\), we have \(g \circ  k \circ  {g}^{-1} \in  P\). And \(P \subseteq  K\), so \(g \circ  k \circ  {g}^{-1} \in  K\), thus
\[
\left( {f \circ  g}\right)  \circ  k \circ  {\left( f \circ  g\right) }^{-1} = f \circ  g \circ  k \circ  {g}^{-1} \circ  {f}^{-1} \in  {E}_{K}\left( P\right).
\]
That is, \(f \circ  g \in  {E}_{K}\left( P\right)\). Therefore, \({E}_{K}\left( P\right)\) is a group under the composition operation of mappings. \hfill $\square$

A special case of this proposition is \(P = K\). At this time, \({E}_{P}\left( P\right)  = N\left( P\right)\). Here, \(N\left( P\right)\) represents the general normalizer of \(P\) in \(\operatorname{Homeo}\left( G\right)\), that is, \(N\left( P\right)  = \{ f \in  \operatorname{Homeo}\left( X\right) : {fP} \subseteq  {Pf}\}\). This result shows that if \(P\) is a set composed of measure-preserving homeomorphisms, then the elements in the normalizer \(N\left( P\right)\) are also measure-preserving homeomorphisms. We will have important results about the general normalizer later. Next, we show that the affine set \(\mathrm{{AF}}\left( G\right)\) can also be represented by the normalizer. That is, when we take \(P = K = \left\{  {{L}_{x} : x \in  G}\right\}\), we can get the affine set \(\operatorname{AF}\left( G\right)\). 

We need the following lemma.

\textbf{Lemma 3.1.1.} Let \(G\) be a compact group and \(\mu\) be the Haar measure on \(G\). Then the normalizer \(N\left( \left\{  {{L}_{x} : x \in  G}\right\}  \right)\) of the left translation \(\left\{  {{L}_{x} : x \in  G}\right\}\) in \(\operatorname{Homeo}\left( G\right)\) is equal to the affine set \(\operatorname{AF}\left( G\right)\) of \(G\).

Proof: For any \(f \in  N\left( \left\{  {{L}_{x} : x \in  G}\right\}  \right)\) and \(y \in  G\), we have
$$
\begin{aligned}
f \circ L_y \circ L_x \circ\left(f \circ L_y\right)^{-1} & =f \circ L_y \circ L_x \circ\left(L_y\right)^{-1} \circ f^{-1} \\
& =f \circ L_y \circ L_x \circ L_{y^{-1}} \circ f^{-1} \\
& =f \circ L_{y x y^{-1}} \circ f^{-1}. 
\end{aligned}
$$
Since \(f \circ  {L}_{{yx}{y}^{-1}} \circ  {f}^{-1} \in  \{ {L}_{x} : x \in  G\}\), so \(f \circ  {L}_{y} \in  N\left( \left\{  {{L}_{x} : x \in  G}\right\}  \right)\).

For any \(f \in  N\left( \left\{  {{L}_{x} : x \in  G}\right\}  \right)\), there exists \(s \in  G\) such that \(f\left( s\right)  = e\), where \(e\) is the identity element, then \(f \circ  {L}_{s}\left( e\right)  = f\left( s\right)  = e\). Let \(g = f \circ  {L}_{s} \in  N\left( \left\{  {{L}_{x} : x \in  G}\right\}  \right)\), then \(g\left( e\right)  = e\) holds for any \(x,z \in  G\). By the definition of \(N\left( \left\{  {{L}_{x} : x \in  G}\right\}  \right)\), there exists an element \(h\left( x\right)  \in  G\) related to \(x\) such that
\[
g \circ  {L}_{x} \circ  {g}^{-1}\left( z\right)  = {L}_{h\left( x\right) }\left( z\right),
\]
Therefore \(g \circ  {L}_{x}\left( z\right)  = {L}_{h\left( x\right) }\left( {g\left( z\right) }\right)\), that is \(g\left( {xz}\right)  = h\left( x\right) g\left( z\right)\). Let \(z = e\), then \(g\left( x\right)  = h\left( x\right) g\left( e\right)  =\)  \(h\left( x\right)\), then \(g\left( {xz}\right)  = g\left( x\right) g\left( z\right)\). Therefore \(g\) is an isomorphism on \(G\). Also, since both the left translation and \(f\) are continuous mappings, \(g = f \circ  {L}_{s}\) is also continuous, that is \(g \in  \operatorname{Aut}\left( G\right)\). Next, we prove that \(f\) is affine. For any \(w \in  G\).
$$
\begin{aligned}
f(w) & =g \circ\left(L_s\right)^{-1}(w)=g \circ L_{s^{-1}}(w) \\
& =g\left(s^{-1} w\right)=g\left(s^{-1}\right) g(w) \\
& =L_{g\left(s^{-1}\right)} \circ g(w),
\end{aligned}
$$
Therefore \(f = {L}_{g\left( {s}^{-1}\right) } \circ  g \in  \mathrm{{AF}}\left( G\right)\).

Conversely, for any \(h \in  \mathrm{{AF}}\left( G\right)\), let \(h = {L}_{y} \circ  g\), where \(y \in  G,g \in  \operatorname{Aut}\left( G\right)\). Then for any \(x \in  G,\)
$$
\begin{aligned}
h \circ L_x \circ h^{-1} & =L_y \circ g \circ L_x \circ\left(L_y \circ g\right)^{-1}=L_y \circ g \circ L_x \circ g^{-1} \circ L_{y^{-1}} \\
& =L_y \circ g \circ L_x \circ L_{g^{-1}\left(y^{-1}\right)} \circ g^{-1}=L_y \circ g \circ L_{x \cdot g^{-1}\left(y^{-1}\right)} \circ g^{-1} \\
& =L_y \circ L_{g(x) \cdot y^{-1}} \circ g \circ g^{-1}=L_{y \cdot g(x) \cdot y^{-1}},
\end{aligned}
$$
Therefore \(h \in  N\left( \left\{  {{L}_{x} : x \in  G}\right\}  \right)\). 

In summary, \(N\left( \left\{  {{L}_{x} : x \in  G}\right\}  \right)  = \operatorname{AF}\left( G\right)\). \hfill $\square$

Since \({E}_{G}\left( \left\{  {{L}_{x} : x \in  G}\right\}  \right)  = N\left( \left\{  {{L}_{x} : x \in  G}\right\}  \right)\), by Lemma 3.1.1, we obtain
\[
{E}_{G}\left( \left\{  {{L}_{x} : x \in  G}\right\}  \right)  = \operatorname{AF}\left( G\right),
\]
Moreover, since we have already proved that \({E}_{G}\left( \left\{  {{L}_{x} : x \in  G}\right\}  \right)\) preserves the Haar measure, this actually gives a proof that affine transformations preserve the Haar measure.

Since the affine transformations of the compact group \(G\) also preserve the Haar measure, we can consider the normalizer \(N\left( {\mathrm{{AF}}\left( G\right) }\right)\) of the affine set and the normalizer \({E}_{G}\left( {\mathrm{{AF}}\left( G\right) }\right)\) of the translation. Since the affine set \(\mathrm{{AF}}\left( G\right)  = {E}_{G}\left( \left\{  {{L}_{x} : x \in  G}\right\}  \right)\), by Proposition 3.1.1, we know that
\[
\mathrm{{AF}}\left( G\right)  = {E}_{G}\left( \left\{  {{L}_{x} : x \in  G}\right\}  \right)  \subseteq  {E}_{G}\left( {\mathrm{{AF}}\left( G\right) }\right) ,
\]
That is, the translational normalizer of an affine set contains the affine set. In fact, the affine set \(\mathrm{{AF}}\left( G\right)\) is also contained in its normalizer \(N\left( {\mathrm{{AF}}\left( G\right) }\right)\), because the affine set \(\mathrm{{AF}}\left( G\right)\) is a subgroup of \(\operatorname{Homeo}\left( G\right)\). The composition of affine transformations and the inverse of an affine transformation are also affine. Combining the above discussion, we actually obtain an inclusion relation between sets composed of mappings that preserve the Haar measure, namely
\[
\left\{  {{L}_{x} : x \in  G}\right\}   \subseteq  \mathrm{{AF}}\left( G\right)  \subseteq  N\left( {\mathrm{{AF}}\left( G\right) }\right)  \subseteq  {E}_{G}\left( {\mathrm{{AF}}\left( G\right) }\right).
\]
We can completely describe this inclusion relation using the generalized normalizer, namely
$$
\begin{aligned}
\left\{L_x: x \in G\right\} & \subseteq E_G\left(\left\{L_x: x \in G\right\}\right) \\
& \subseteq E_{E_{G}\left\{L_x: x \in G\right\}}\left(E_{G}\left\{L_x: x \in G\right\}\right) \\
& \subseteq E_G\left(E_{G}\left\{L_x: x \in G\right\}\right).
\end{aligned}
$$
From this, we can see that the generalized normalizer actually covers many known results, so it is a very valuable concept. If we can show that the affine set \(\mathrm{{AF}}\left( G\right)\) of some compact group \(G\) is properly contained in \({E}_{G}\left( {\mathrm{{AF}}\left( G\right) }\right)\), it means that we can indeed find new self-homeomorphisms that preserve the Haar measure using this concept.

Next, we use a constructive method to show that in some non-commutative compact groups, the affine set \(\mathrm{{AF}}\left( G\right)\) is indeed properly contained in \({E}_{G}\left( {\mathrm{{AF}}\left( G\right) }\right)\). Before proving this result, we first need to introduce some properties of the translational normalizer \({E}_{G}\left( {\operatorname{AF}\left( G\right) }\right)\) of the affine set itself.

\textbf{Lemma 3.1.2.} Let \(f\) be a self - homeomorphism on \(G\) and \(f\left( e\right)  = e\), then \(f \in  {E}_{G}\left( {\operatorname{AF}\left( G\right) }\right)\) if and only if for any \(x \in  G\), there exists \({\varphi }_{x} \in  \operatorname{Aut}\left( G\right)\) such that for any \(z \in  G\), \(f\left( {xz}\right)  = f\left( x\right) {\varphi }_{x}\left( {f\left( z\right) }\right)\) holds. At this time, there exists a group homomorphism \(\phi\) from \(G\) to \(\operatorname{Aut}\left( G\right)\) such that for any \(x \in  G\), \(\phi \left( x\right)  = {\varphi }_{x}\) holds.

At this time, for any \(x,y \in  G\), the self-homeomorphism \({\varphi }_{x},{\varphi }_{y}\) related to \(x,y\) satisfies \({\varphi }_{x} \circ  {\varphi }_{y} = {\varphi }_{xy}\), and \({\varphi }_{e}\) is the identity mapping.

Proof: Necessity: For any \(f \in  {E}_{G}\left( {\mathrm{{AF}}\left( G\right) }\right)\) and \(f\left( e\right)  = e\) and \(x \in  G\), there exist left translations \({L}_{x}\) and \({\varphi }_{x} \in  \operatorname{Aut}\left( G\right)\) such that \(f \circ  {L}_{x} \circ  {f}^{-1} = {L}_{y} \circ  {\varphi }_{x}\). For any \(z \in  G\),
$$
\begin{aligned}
f \circ L_x \circ f^{-1}(z) & =L_y \circ \varphi_x(z) \\
\Rightarrow  f \circ L_x(z) & =L_y \circ \varphi_x(f(z)) \\
\Rightarrow  f(x z) & =y \varphi_x(f(z)),
\end{aligned}
$$
Take \(z = e\), then \(f\left( x\right)  = y{\varphi }_{x}\left( {f\left( e\right) }\right)  = y{\varphi }_{x}\left( e\right)  = y\), so \(f\left( {xz}\right)  = f\left( x\right) {\varphi }_{x}\left( {f\left( z\right) }\right)\). Sufficiency: For any \(x \in  G\), there exists \({\varphi }_{x} \in  \operatorname{Aut}\left( G\right)\) such that for any \(z \in  G\), there is \(f\left( {xz}\right)  = f\left( x\right) {\varphi }_{x}\left( {f\left( z\right) }\right)\). Since \(f\) is a bijection, then for any \(z \in  G,f\left( {x{f}^{-1}\left( z\right) }\right)  = f\left( x\right) {\varphi }_{x}\left( z\right)\), that is \(f \circ  {L}_{x} \circ  {f}^{-1} =\)  \({L}_{f\left( x\right) } \circ  {\varphi }_{x}\), namely \(f \in  E\left( G\right)\).

For any \(x,y,z \in  G\), there is
$$
\begin{aligned}
f(x y z) & =f(x) \varphi_x(f(y z)) \\
& =f(x) \varphi_x\left(f(y) \varphi_y(f(z))\right) \\
& =f(x) \varphi_x(f(y)) \varphi_x \circ \varphi_y(f(z)),
\end{aligned}
$$
On the other hand, there is
$$
\begin{aligned}
f(x y z) & =f(x y) \varphi_{x y}(f(z)) \\
& =f(x) \varphi_x(f(y)) \varphi_{x y}(f(z)).
\end{aligned}
$$
By comparison, for any \(x,y,z \in  G\), there is \({\varphi }_{x} \circ  {\varphi }_{y}\left( {f\left( z\right) }\right)  = {\varphi }_{xy}\left( {f\left( z\right) }\right)\). Since \(f\) is a bijection, so for any \(x,y,z \in  G,{\varphi }_{x} \circ  {\varphi }_{y}\left( z\right)  = {\varphi }_{xy}\left( z\right)\), that is \({\varphi }_{x} \circ  {\varphi }_{y} = {\varphi }_{xy}\). Let \(x = y = e\), then \({\varphi }_{e} \circ  {\varphi }_{e} = {\varphi }_{e}\), so \({\varphi }_{e}\) is the identity mapping. This shows that there exists a group homomorphism \(\phi\) from \(G\) to \(\mathrm{{AF}}\left( G\right)\) such that for any \(x \in  G\), there is \(\phi \left( x\right)  = {\varphi }_{x}\). \hfill $\square$

Next, we will show that in some non-commutative compact groups, the translation normalizer \({E}_{G}\left( {\operatorname{AF}\left( G\right) }\right)\) of the affine set can indeed strictly expand the affine set \(\mathrm{{AF}}\left( G\right)\). The following theorem gives a specific method for constructing the mapping in \({E}_{G}\left( {\mathrm{{AF}}\left( G\right) }\right)\).

\textbf{Theorem 3.1.2.} Let \(g\) and \(h\) be continuous endomorphisms on a non-commutative compact group \(G\). For any \(z \in  G\), let \(f\left( z\right)  = g\left( z\right) h{\left( z\right) }^{-1}\). If \(f\) is a homeomorphism on \(G\), then \(f \in  {E}_{G}\left( {\operatorname{AF}\left( G\right) }\right)\).

Proof: For any \(z \in  G\), let \({\varphi }_{x}\left( z\right)  = h\left( x\right) {zh}{\left( x\right) }^{-1}\). Since \(h\) is continuous, by the closure of the topological group under group operations and inverse operations, \({\varphi }_{x}\) is continuous. Also, since \({\varphi }_{x}\left( {h{\left( x\right) }^{-1}{zh}\left( x\right) }\right)  = z\), \({\varphi }_{x}\) is a continuous
bijection, so \({\varphi }_{x}\) is a homeomorphism on \(G\). For any \(z,w \in  G\),
$$
\begin{aligned}
\varphi_x(z w) & =h(x) z w h(x)^{-1} \\
& =h(x) z h(x)^{-1} h(x) w h(x)^{-1} \\
& =\varphi_x(z) \varphi_x(w).
\end{aligned}
$$
Therefore, \({\varphi }_{x} \in  \operatorname{Aut}\left( G\right)\). For any \(z \in  G\), since \(g\) and \(h\) are endomorphisms on \(G\), then for any \(x,z \in  G\),
$$
\begin{aligned}
f(x z) & =g(x z) h(x z)^{-1} \\
& =g(x) g(z)(h(x) h(z))^{-1} \\
& =g(x) g(z) h(z)^{-1} h(x)^{-1} \\
& =g(x) h(x)^{-1} h(x) g(z) h(z)^{-1} h(x)^{-1} \\
& =f(x) h(x) f(z) h(x)^{-1} \\
& =f(x) \varphi_x(f(z)).
\end{aligned}
$$
Finally, \(f\left( e\right)  = g\left( e\right) h{\left( e\right) }^{-1} = e\). By Lemma 3.1.2, \(f \in  {E}_{G}\left( {\mathrm{{AF}}\left( G\right) }\right)\). \hfill $\square$

Since the proof of Theorem 3.1.2 is constructive, we can use this theorem to construct specific examples. We will use a specific example below to show that in some non-commutative compact groups, this method can yield homeomorphisms that preserve the Haar measure other than affine ones. We consider the direct product group \({SO}\left( 3\right)  \times  {SO}\left( 3\right)\) of the special orthogonal group of order three with itself.

\textbf{Example 3.1.1.} For any \(\left( {A,B}\right)  \in  {SO}\left( 3\right)  \times  {SO}\left( 3\right)\), let \(f\left( {A,B}\right)  = \left( {A{B}^{T},B}\right)\), then \(f\) is a non-affine measure-preserving homeomorphism.

Proof: Denote the identity matrix by \(E\). For any \(\left( {A,B}\right)  \in  {SO}\left( 3\right)  \times  {SO}\left( 3\right)\), let \(g\left( {A,B}\right)  =\)  \(\left( {A,B}\right) ,h\left( {A,B}\right)  = \left( {B,E}\right)\), then \(f\left( {A,B}\right)  = \left( {A,B}\right) \left( {{B}^{T},E}\right)  = g\left( {A,B}\right) h{\left( A,B\right) }^{-1}\). For any \(\left( {{A}_{1},{B}_{1}}\right)\) and \(\left( {{A}_{2},{B}_{2}}\right)\),
$$
\begin{aligned}
g\left(\left(A_1, B_1\right)\left(A_2, B_2\right)\right) & =g\left(A_1 A_2, B_1 B_2\right) \\
& =\left(A_1 A_2, B_1 B_2\right) \\
& =\left(A_1, A_2\right)\left(B_1, B_2\right) \\
& =g\left(A_1, B_1\right) g\left(A_2, B_2\right),
\end{aligned}
$$
On the other hand,
$$
\begin{aligned}
h\left(\left(A_1, B_1\right)\left(A_2, B_2\right)\right) & =h\left(A_1 A_2, B_1 B_2\right) \\
& =\left(B_1 B_2, E\right) \\
& =\left(B_1, E\right)\left(B_2, E\right) \\
& =h\left(A_1, B_1\right) h\left(A_2, B_2\right).
\end{aligned}
$$
Therefore, \(g\) and \(h\) are endomorphisms on \({SO}\left( 3\right)  \times  {SO}\left( 3\right)\). For any \(\left( {A,B}\right) ,f\left( {{AB},B}\right)  = \left( {{AB}{B}^{T},B}\right)  =\) (A, B), since \(f\) is a bijection, and because the topological group is continuous with respect to group operations and inverse operations, \(f\) is continuous. Similarly, \({f}^{-1}\) is continuous. Thus, \(f\) is a homeomorphism on \({SO}\left( 3\right)  \times  {SO}\left( 3\right)\). By Theorem 3.1.2, \(f \in  {E}_{G}(\mathrm{{AF}}({SO}\left( 3\right)  \times\)  \({SO}\left( 3\right) ))\), so \(f\) preserves the Haar measure. Next, we prove \(f \notin  \mathrm{{AF}}\left( {{SO}\left( 3\right)  \times  {SO}\left( 3\right) }\right)\) by contradiction. Assume \(f \in  \operatorname{AF}\left( {{SO}\left( 3\right)  \times  {SO}\left( 3\right) }\right)\). Since \(f\left( {E,E}\right)  = \left( {E,E}\right)\), then \(f \in  \operatorname{Aut}\left( G\right)\). So, for any \({B}_{1},{B}_{2} \in  {SO}\left( 3\right)\),
$$
\begin{aligned}
f\left(\left(E, B_1\right)\left(E, B_2\right)\right) & =f\left(E, B_1 B_2\right) \\
& =\left(\left(B_1 B_2\right)^T, B_1 B_2\right) \\
& =\left(B_2{ }^T B_1^T, B_1 B_2\right),
\end{aligned}
$$
On the other hand,
$$
\begin{aligned}
f\left(\left(E, B_1\right)\left(E, B_2\right)\right) & =f\left(E, B_1\right) f\left(E, B_2\right) \\
& =\left(B_1^T, B_1\right)\left(B_2^T, B_2\right) \\
& =\left(B_1^T B_2{ }^T, B_1 B_2\right).
\end{aligned}
$$
Therefore, \({B}_{2}{}^{T}{B}_{1}^{T} = {B}_{1}{}^{T}{B}_{2}{}^{T}\). Taking the transpose of both sides, we get \({B}_{2}{B}_{1} = {B}_{1}{B}_{2}\), which contradicts the fact that \({SO}\left( 3\right)\) is a non - abelian group! Thus, \(f \notin  \mathrm{{AF}}\left( {{SO}\left( 3\right)  \times  {SO}\left( 3\right) }\right)\). \hfill $\square$

From this example, we can see that for some non-abelian compact groups, the translation normalizer \({E}_{G}\left( {\mathrm{{AF}}\left( G\right) }\right)\) can indeed expand the affine set \(\mathrm{{AF}}\left( G\right)\), which also illustrates the significance of this concept. Due to the monotonicity of the normalizer, we can continuously consider the normalizer of the normalizer, which may continuously expand the affine set, forming a larger set of mappings that contains the affine set. This is beneficial to analyzing the structure of the homeomorphism group \({\operatorname{Homeo}}_{\mu }\left( G\right)\) that preserves the Haar measure to some extent.

\subsection{The Relationship between the Affine Set and its Normalizer}\

In the previous section, we showed that for some non-commutative compact groups, the translation normalizer of their affine sets strictly contains the affine sets. However, not all compact groups are like this. But the fact that the translation normalizer is exactly equal to the affine set is not meaningless. On the contrary, we can obtain a certain characterization of the structure of the group of measure-preserving homeomorphisms from this, which enables us to have a deeper understanding of the group of self-homeomorphisms \({\operatorname{Homeo}}_{\mu }\left( G\right)\) that preserve the Haar measure. In this section, we consider some commutative compact groups and prove that the translation normalizer of the affine sets of these commutative compact groups is actually the affine set.

First, consider the finite cyclic group \(\mathbb{Z}/m\mathbb{Z}\) with the discrete topology. Since any mapping is continuous under the discrete topology, the isomorphism of \(\mathbb{Z}/m\mathbb{Z}\) is just the general group isomorphism, and the affine is the affine of the group.

\textbf{Proposition 3.2.1.} \(\operatorname{AF}\left( {\mathbb{Z}/m\mathbb{Z}}\right)  = N\left( {\operatorname{AF}\left( {\mathbb{Z}/m\mathbb{Z}}\right) }\right)  = {E}_{\mathbb{Z}/m\mathbb{Z}}\left( {\operatorname{AF}\left( {\mathbb{Z}/m\mathbb{Z}}\right) }\right)\).

Proof: For any \(f \in  {E}_{\mathbb{Z}/m\mathbb{Z}}\left( {\operatorname{AF}\left( {\mathbb{Z}/m\mathbb{Z}}\right) }\right)\) such that \(f\left( 0\right)  = 0\), then for any \(x,z \in\)  \(\mathbb{Z}/m\mathbb{Z}\), there exists an automorphism \({\phi }_{x}\) such that \(f\left( {xz}\right)  = f\left( x\right)  + {\phi }_{x}\left( {f\left( z\right) }\right)\). Take \(x = 1\), then \(f\left( z\right)  = f\left( 1\right)  +\)  \({\phi }_{1}\left( {f\left( {z - 1}\right) }\right)\). Since \(\operatorname{Aut}\left( {\mathbb{Z}/m\mathbb{Z}}\right)\) is isomorphic to the multiplicative group \({\left( \mathbb{Z}/m\mathbb{Z}\right) }^{ * }\), then \({\phi }_{1}\left( z\right)  = {tz}\), where
\(t \in  {\left( \mathbb{Z}/m\mathbb{Z}\right) }^{ * }\), let \(f\left( 1\right)  = s\), then \(f\left( {1 + z}\right)  = s + {tf}\left( z\right)\). Since \(f\left( 0\right)  = 0\), then
$$
\begin{aligned} 
f(z) & =s+t f(z-1) \\
& =s+t(s+t f(z-2)) \\
& =(1+t) s+t^2 f(z-2) \\
& =\left(1+t+t^2\right) s+t^3 f(z-3) \\
& =\cdots \cdots \cdots \\
& =\left(1+t+t^2+\cdots+t^{z-1}\right) s+t^z f(0) \\
& =\left(1+t+t^2+\cdots+t^{z-1}\right) s.
\end{aligned}
$$
Next, we will prove \(t = 1\) by contradiction. If \(t \neq  1\), since \(1 + t + {t}^{2} + \cdots  + {t}^{z - 1} = \frac{{t}^{2} - 1}{t - 1}\) and \(f\left( 0\right) ,f\left( 1\right) ,f\left( 2\right) ,\cdots ,f\left( {n - 1}\right)\) are pairwise distinct, this is equivalent to \(0,t - 1,{t}^{2} - 1,\cdots ,{t}^{n - 1} - 1\) being pairwise distinct, which is also equivalent to \(1,t,{t}^{2},\cdots ,{t}^{n - 1}\) being pairwise distinct. And since \(t \in  {\left( \mathbb{Z}/m\mathbb{Z}\right) }^{ * }\), then \(1,t,{t}^{2},\cdots ,{t}^{n - 1}\) are all elements of \({\left( \mathbb{Z}/m\mathbb{Z}\right) }^{ * }\). Therefore, there are at most Euler's totient number $\phi \left( n\right)$ choices, which contradicts the fact that \(1,t,{t}^{2},\cdots ,{t}^{n - 1}\) are pairwise distinct. Thus, \(t = 1\) , then \(f\left( z\right)  = {zs}\). Since \(f\left( {z}_{1}\right)  =\)  \(f\left( {z}_{2}\right)  \Leftrightarrow  \left( {{z}_{1} - {z}_{2}}\right) s = 0\left( {\;\operatorname{mod}\;n}\right)\), \(f\) is an injective function if and only if \(n\) and \(s\) are co-prime. In this case, \(f\) is a bijective function and \(t \in  \left( {\mathbb{Z}/m\mathbb{Z}}\right) {)}^{ * }\), that is \(f \in  \operatorname{Aut}\left( {\mathbb{Z}/m\mathbb{Z}}\right)\). Since for any \(g \in  {E}_{\mathbb{Z}/m\mathbb{Z}}\left( {\operatorname{AF}\left( {\mathbb{Z}/m\mathbb{Z}}\right) }\right)\) , there exists a translation \({L}_{g{\left( 0\right) }^{-1}}\) such that \({L}_{g{\left( 0\right) }^{-1}} \circ  g \in  {E}_{\mathbb{Z}/m\mathbb{Z}}\left( {\operatorname{AF}\left( {\mathbb{Z}/m\mathbb{Z}}\right) }\right)\) and \({L}_{g{\left( 0\right) }^{-1}} \circ  g\left( 0\right)  = 0\). Therefore, \({L}_{g{\left( 0\right) }^{-1}} \circ  g \in  \operatorname{Aut}\left( {\mathbb{Z}/m\mathbb{Z}}\right)\), then \(g \in  \operatorname{AF}\left( {\mathbb{Z}/m\mathbb{Z}}\right)\).

In conclusion, \({E}_{\mathbb{Z}/m\mathbb{Z}}\left( {\mathrm{{AF}}\left( {\mathbb{Z}/m\mathbb{Z}}\right) }\right)  \subseteq  \mathrm{{AF}}\left( {\mathbb{Z}/m\mathbb{Z}}\right)\) , and in Section 3.1, we have already explained
\[
\operatorname{AF}\left( G\right)  \subseteq  N\left( {\operatorname{AF}\left( G\right) }\right)  \subseteq  {E}_{G}\left( {\operatorname{AF}\left( G\right) }\right),
\]
Therefore \(\operatorname{AF}\left( {\mathbb{Z}/m\mathbb{Z}}\right)  = N\left( {\operatorname{AF}\left( {\mathbb{Z}/m\mathbb{Z}}\right) }\right)  = {E}_{\mathbb{Z}/m\mathbb{Z}}\left( {\operatorname{AF}\left( {\mathbb{Z}/m\mathbb{Z}}\right) }\right)\). \hfill $\square$    

Next, we consider the \(n\)-dimensional torus \({\mathbb{T}}^{n}\).

\textbf{Theorem 3.2.1.} \(\operatorname{AF}\left( {\mathbb{T}}^{n}\right)  = N\left( {\operatorname{AF}\left( {\mathbb{T}}^{n}\right) }\right)  = {E}_{{\mathbb{T}}^{n}}\left( {\operatorname{AF}\left( {\mathbb{T}}^{n}\right) }\right)\).

Proof: Let \(f \in  {E}_{{\mathbb{T}}^{n}}\left( {\mathrm{{AF}}\left( {\mathbb{T}}^{n}\right) }\right)\), and \(f\left( {1,1,\cdots ,1}\right)  = \left( {1,1,\cdots ,1}\right)\). For any \(x =\)  \(\left( {{x}_{1},{x}_{2},\cdots ,{x}_{n}}\right)  \in  {\mathbb{T}}^{n}\), by the definition of \({E}_{{\mathbb{T}}^{n}}\left( {\mathrm{{AF}}\left( {\mathbb{T}}^{n}\right) }\right)\), there exist a left translation \({L}_{y\left( x\right) }\) related to \(x\) and \({\varphi }_{x} \in  \operatorname{Aut}\left( {\mathbb{T}}^{n}\right)\) such that for any \(z = \left( {{z}_{1},{z}_{2},\cdots ,{z}_{n}}\right)  \in  {\mathbb{T}}^{n}\), we have
\[
f \circ  {L}_{x} \circ  {f}^{-1}\left( {{z}_{1},{z}_{2},\cdots ,{z}_{n}}\right)  = y\left( x\right) {\varphi }_{x}\left( {{z}_{1},{z}_{2},\cdots ,{z}_{n}}\right),
\]
Let \(z = \left( {{z}_{1},{z}_{2},\cdots ,{z}_{n}}\right)  = \left( {1,1,\cdots ,1}\right)\) , then we get \(y\left( x\right)  = f\left( x\right)\) . Therefore
\[
f \circ  {L}_{x} \circ  {f}^{-1}\left( {{z}_{1},{z}_{2},\cdots ,{z}_{n}}\right)  = f\left( x\right) {\varphi }_{x}\left( {{z}_{1},{z}_{2},\cdots ,{z}_{n}}\right),
\]
That is
\[
f{\left( x\right) }^{-1} \cdot  f \circ  {L}_{x} \circ  {f}^{-1}\left( {{z}_{1},{z}_{2},\cdots ,{z}_{n}}\right)  = {\varphi }_{x}\left( {{z}_{1},{z}_{2},\cdots ,{z}_{n}}\right).
\]
According to Lemma 3.1.2,
\[
{\varphi }_{x} = {\varphi }_{\left( {x}_{1},{x}_{2},\cdots ,{x}_{n}\right) } = {\varphi }_{\left( {x}_{1},1,\cdots ,1\right) } \circ  {\varphi }_{\left( 1,{x}_{2},\cdots ,1\right) } \circ  {\varphi }_{\left( 1,1,\cdots ,{x}_{n}\right) },
\]
Therefore, we only need to consider the case where only one component of \(x\) is not equal to 1. For any \(l \in  \{ 1,2,\cdots ,n\}\), let \({t}_{l} = \left( {1,\cdots ,1,{x}_{l},1,\cdots ,1}\right)\). Because the automorphisms of a torus correspond one-to-one with matrices., there exists
integers 
$${k}_{l11}\left( {x}_{l}\right) ,{k}_{l12}\left( {x}_{l}\right) ,\cdots ,{k}_{nn}\left( {x}_{l}\right) ,\cdots ,{k}_{ln1}\left( {x}_{l}\right) ,\cdots ,{k}_{lnn}\left( {x}_{l}\right)$$
such that
\[
{\varphi }_{{t}_{l}}\left( {{z}_{1},\cdots ,{z}_{n}}\right)  = \left( {{z}_{1}^{{k}_{l11}\left( {x}_{l}\right) }\cdots {z}_{n}^{{k}_{l1n}\left( {x}_{l}\right) },\cdots ,{z}_{1}^{{k}_{ln1}\left( {x}_{l}\right) }\cdots {z}_{n}^{{k}_{lnn}\left( {x}_{l}\right) }}\right),
\]
Take \({z}_{1} = \cdots {z}_{j - 1} = {z}_{j + 1} = \cdots {z}_{n} = 1\), then
\[
f{\left( {t}_{l}\right) }^{-1} \cdot  f \circ  {L}_{{t}_{l}} \circ  {f}^{-1}\left( {1,\cdots ,1,{z}_{j},1,\cdots ,1}\right)  = \left( {{z}_{j}^{{k}_{lj1}\left( {x}_{l}\right) },\cdots ,{z}_{j}^{{k}_{ljn}\left( {x}_{l}\right) }}\right).
\]

Since \(f\) is continuous with respect to \({t}_{l} \in  {\mathbb{T}}^{n}\), it is continuous with respect to \({x}_{l} \in  \mathbb{T}\). For any \({z}_{j} \in  \mathbb{T}\), by the continuity of the topological group, \(f{\left( {x}_{l}\right) }^{-1} \cdot  f \circ  {L}_{{x}_{l}} \circ  {f}^{-1}\left( {1,\cdots ,1,{z}_{j},1,\cdots ,1}\right)\) is continuous with respect to \({x}_{l}\) , then \({z}_{j}^{{k}_{lj1}\left( {x}_{l}\right) },{z}_{j}^{{k}_{lj2}\left( {x}_{l}\right) },\cdots ,{z}_{j}^{{k}_{ljn}\left( {x}_{l}\right) }\) is continuous with respect to \({x}_{l}\) respectively. Next, we prove that \({\left\{  {k}_{lji}\left( {x}_{l}\right) \right\}  }_{i,j,l = 1}^{n}\) is a constant function on \({X}_{l} \in  \mathbb{T}\). We use the method of proof by contradiction. Assume that there exists \(i,j,l \in  \{ 1,2,\cdots ,n\}\) such that \({k}_{lji}\left( {x}_{l}\right)\) is not a constant function on \(\mathbb{T}\). First, we prove that there exists \({x}_{l0} \in  \mathbb{T}\) such that \({k}_{lji}\left( {x}_{l}\right)\) is not a constant in any neighborhood of \({x}_{{l}_{0}}\). Otherwise, for any \({x}_{l} \in  \mathbb{T}\), there exists a neighborhood \({A}_{{x}_{l}}\) of \({x}_{l}\) such that \({k}_{lji}\left( {x}_{l}\right)\) is a constant on \({A}_{{x}_{l}}\). Take an open interval \({U}_{{x}_{l}}\) with arc length \({\delta }_{x}\left( {{\delta }_{x} \leq  1/2}\right)\) on the unit circle centered at \({x}_{l}\) such that \({U}_{{x}_{l}} \subseteq  {A}_{{x}_{l}}\). Since \(\mathbb{T}\) is a compact set, there exists \({x}_{l1},{x}_{l2},\cdots ,{x}_{ls} \in  \mathbb{T}\) such that \({U}_{{x}_{l1}} \cup  {U}_{{x}_{l2}} \cup  \cdots  \cup  {U}_{{x}_{ls}} = \mathbb{T}\). Since \(\mathbb{T}\) is a connected set, for \({U}_{{x}_{l1}}\), there must exist \({U}_{{x}_{lt}}\) such that \({U}_{{x}_{l1}} \cap  {U}_{{x}_{lt}}\) is non-empty. Therefore, \({k}_{lji}\left( {x}_{j}\right)\) is a constant on \({U}_{{x}_{l1}} \cup  {U}_{{x}_{lt}}\). For \({U}_{{x}_{l1}} \cup  {U}_{{x}_{lt}}\), there must exist \({U}_{{x}_{lu}}\) such that \(\left( {{U}_{{x}_{l1}} \cup  {U}_{{x}_{lt}}}\right)  \cap  {U}_{{x}_{lu}}\) is non-empty. Therefore, \({k}_{jji}\left( {x}_{l}\right)\) is a constant on \({U}_{{x}_{l1}} \cup  {U}_{{x}_{{k}_{lt}}} \cup  {U}_{{x}_{lu}}\). Repeating this process, since \({U}_{{x}_{l1}},{U}_{{x}_{l2}},\cdots ,{U}_{{x}_{ls}}\) is finite, it will stop after a certain step. Then \({k}_{lji}\left( {x}_{l}\right)\) is a constant on \(\mathbb{T}\), which contradicts the fact that \({k}_{lji}\left( {x}_{l}\right)\) is not a constant function on \(\mathbb{T}\) . Therefore, there must exist \({x}_{{l}_{0}} \in  \mathbb{T}\) such that \({k}_{lji}\left( {x}_{l}\right)\) is not a constant in any neighborhood of \({x}_{{l}_{0}}\). Since \(i,j,l \in  \{ 1,2,\cdots ,n\} ,{z}_{j}^{{k}_{lji}\left( {x}_{l}\right) }\) is continuous with respect to \({x}_{l} \in  \mathbb{T}\) for any \(i,j,l \in  \{ 1,2,\cdots ,n\} ,{z}_{j}^{{k}_{lji}\left( {x}_{l}\right) }\), when \({x}_{l} \rightarrow  {x}_{{l}_{0}}\), \({z}_{j}^{{k}_{lji}\left( {x}_{l}\right) } \rightarrow  {z}_{j}^{{k}_{lji}\left( {x}_{l0}\right) }\), that is, when \({x}_{l} \rightarrow  {x}_{{l}_{0}}\) , \({z}_{j}^{{k}_{lji}\left( {x}_{l}\right)  - {k}_{lji}\left( {x}_{{l}_{0}}\right) } \rightarrow  1\). Since \({k}_{lji}\left( {x}_{l}\right)\) is not a constant in any neighborhood of \({x}_{{l}_{0}}\), there exists an integer sequence \({\left\{  {k}_{n}\right\}  }_{n = 1}^{\infty }\) with each term non-zero such that \({z}^{{k}_{n}} \rightarrow  1\left( {n \rightarrow  \infty }\right)\) for any \(z \in  \mathbb{T}\).

We will now prove \(\mu \left( \left\{  {z \in  \mathbb{T} : {z}^{{k}_{n}} \rightarrow  1,n \rightarrow  \infty }\right\}  \right)  = 0\) to derive a contradiction. Consider the set
$$
\left\{z \in \mathbb{T}: z^{k_n} \rightarrow 1, n \rightarrow \infty\right\}=\bigcap_{s=2}^{\infty} \bigcup_{N=1}^{\infty} \bigcap_{n=N}^{\infty}\left\{z \in \mathbb{T}: d\left(z^{k_n}-1\right) \leq \frac{1}{s}\right\}
$$
For any \(n \geq  N\),
$$
\begin{aligned}
& \left\{z \in \mathbb{T}: d\left(z^{k_n}-1\right) \leq 1 / s\right\}=\left\{e^x: x \in[0,1), d\left(e^{x k_n}-1\right) \leq \frac{1}{s}\right\} \\
& =\left\{e^x: x \in\left[0, \frac{1}{s k_n}\right) \cup\left[\frac{1}{k_n}, \frac{1+s}{s k_n}\right) \cup \cdots \cup\left[\frac{k_n-1}{k_n}, \frac{1+\left(k_n-1\right) s}{s k_n}\right) \cup\right. \\
& \left.\quad\left(\frac{s-1}{s k_n}, \frac{1}{k_n}\right) \cup\left(\frac{2 s-1}{s k_n}, \frac{2}{k_n}\right) \cup \cdots \cup\left(\frac{k_n s-1}{s k_n}, 1\right)\right\},
\end{aligned}
$$
Then \(\mu \left( \left\{  {z \in  \mathbb{T} : d\left( {{z}^{{k}_{n}} - 1}\right)  \leq  \frac{2}{s}}\right\}  \right)  = \frac{1}{s{k}_{n}} \times  {k}_{n} + \frac{1}{s{k}_{n}} \times  {k}_{n} = \frac{2}{s}\), so
$$
\begin{aligned}
\mu\left(\left\{z \in \mathbb{T}: z^{k_n} \rightarrow 1\right\}\right) & =\mu\left(\bigcap_{s=2}^{\infty} \bigcup_{N=1}^{\infty} \bigcap_{n=N}^{\infty}\left\{z \in \mathbb{T}: d\left(z^{k_n}-1\right) \leq \frac{1}{s}\right\}\right) \\
& \leq \mu\left(\bigcap_{s=2}^{\infty} \bigcup_{N=1}^{\infty}\left\{z \in \mathbb{T}: d\left(z^{k_N}-1\right) \leq \frac{1}{s}\right\}\right) \\
& =\mu\left(\bigcap_{s=2}^{\infty}\left\{z \in \mathbb{T}: d\left(z^{k_1}-1\right) \leq \frac{1}{s}\right\}\right) \\
& =\lim _{s \rightarrow \infty} \mu\left(\left\{z \in \mathbb{T}: d\left(z^{k_1}-1\right) \leq \frac{1}{s}\right\}\right) \\
& =\lim _{s \rightarrow \infty} \frac{2}{s}=0.
\end{aligned}
$$
However, for any \(z \in  \mathbb{T}\), there is always \({z}^{{k}_{n}} \rightarrow  1\left( {n \rightarrow  \infty }\right)\). Therefore, \(\left\{  {z \in  \mathbb{T} : {z}^{{k}_{n}} \rightarrow  1,n \rightarrow  \infty }\right\}   = \mathbb{T}\), and then \(\mu \left( \left\{  {z \in  \mathbb{T} : {z}^{{k}_{n}} \rightarrow  1,n \rightarrow  \infty }\right\}  \right)  = 1\), which is a contradiction! Thus, for any \(i,j,l \in  \{ 1,2,\cdots ,n\}\), \({k}_{lji}\left( {x}_{l}\right)\) is a constant function on \({x}_{l} \in  \mathbb{T}\). Denote it as \({k}_{lji}\left( {x}_{l}\right)  = {k}_{lji}\), that is
$$
\begin{aligned}
\varphi_{\left(1, \cdots, 1, x_l, 1, \cdots, 1\right)}\left(z_1, \cdots, z_n\right) & =\varphi_{t_l}\left(z_1, z_2, \cdots, z_n\right) \\
& =\left(z_1^{k_{l 11}} z_2^{k_{112}} \cdots z_n^{k_{l 1 n}}, \cdots, z_1^{k_{l n}} z_2^{k_{l n 2}} \cdots z_n^{k_{l n n}}\right).
\end{aligned}
$$
Moreover, from \({\varphi }_{x} = {\varphi }_{\left( {x}_{1},{x}_{2},\cdots ,{x}_{n}\right) } = {\varphi }_{\left( {x}_{1},1,\cdots ,1\right) } \circ  {\varphi }_{\left( 1,{x}_{2},\cdots ,1\right) } \circ  {\varphi }_{\left( 1,1,\cdots ,{x}_{n}\right) }\), there exists an integer \({k}_{11},\cdots ,{k}_{nn}\) such that for any \(x = \left( {{x}_{1},{x}_{2},\cdots ,{x}_{n}}\right)  \in  {\mathbb{T}}^{n}\) and any \(\left( {{z}_{1},{z}_{2},\cdots ,{z}_{n}}\right)  \in  {\mathbb{T}}^{n}\), there is
$$
\varphi_x\left(z_1, \cdots, z_n\right)=\varphi_x\left(z_1, \cdots, z_n\right)=\left(z_1^{k_{11}}  \cdots z_n^{k_{10}}, \cdots, z_1^{k_{11}}  \cdots z_n^{k_{n n}}\right).
$$
According to Lemma 3.1.2, \({\varphi }_{\left( 1,1,\cdots ,1\right) }\) is the identity mapping. Therefore, for any \(\left( {{z}_{1},{z}_{2},\cdots ,{z}_{n}}\right)  \in  {\mathbb{T}}^{n}\),
\(\left( {{z}_{1},\cdots ,{z}_{n}}\right)  = {\varphi }_{\left( 1,1,\cdots ,1\right) }\left( {{z}_{1},\cdots ,{z}_{n}}\right)  = \left( {{z}_{1}^{{k}_{11}}\cdots {z}_{n}^{{k}_{1n}},\cdots ,{z}_{1}^{{k}_{n1}}\cdots {z}_{n}^{{k}_{nn}}}\right),\)

Take \(\left( {{z}_{1},\cdots ,{z}_{n}}\right)  = \left( {1,\cdots ,1,i,1,\cdots 1}\right)\) , where the \(j\)-th component is \(i\) and the remaining components are 1. Then
\[
\left( {1,\cdots ,1,i,1,\cdots 1}\right)  = \left( {{1}^{{k}_{j1}},\cdots ,{i}^{jj},\cdots ,{1}^{{k}_{jn}}}\right),
\]
From this, we obtain
\[
{k}_{ji} = \left\{  \begin{array}{ll} 1 & i = j \\  0 & i \neq  j \end{array}\right.
\]
That is, for any \(x = \left( {{x}_{1},{x}_{2},\cdots ,{x}_{n}}\right)  \in  {\mathbb{T}}^{n},{\varphi }_{x}\), it is the identity mapping. Therefore, for any \(x \in  {\mathbb{T}}^{n}\) and \(z \in  {\mathbb{T}}^{n}\),
\[
f\left( {\left( {{x}_{1},{x}_{2},\cdots ,{x}_{n}}\right)  \cdot  \left( {{z}_{1},{z}_{2},\cdots ,{z}_{n}}\right) }\right)  = f\left( {{x}_{1},{x}_{2},\cdots ,{x}_{n}}\right) f\left( {{z}_{1},{z}_{2},\cdots ,{z}_{n}}\right).
\]
That is \(f \in  \operatorname{Aut}\left( {\mathbb{T}}^{n}\right)\). For any \(g \in  \operatorname{Homeo}\left( {\mathbb{T}}^{n}\right)\), if \(g\left( {1,1,\cdots ,1}\right)  \neq  \left( {1,1,\cdots ,1}\right)\), then let \(f = {L}_{g{\left( 1,1,\cdots ,1\right) }^{-1}} \circ  g\), then \(f\left( {1,1,\cdots ,1}\right)  = \left( {1,1,\cdots ,1}\right)\). From the above discussion \(f \in  \operatorname{Aut}\left( {\mathbb{T}}^{n}\right)\), thus \(g = {L}_{g\left( {1,1,\cdots ,1}\right) } \circ  f \in  \mathrm{{AF}}\left( G\right)\).

In summary, \({E}_{{\mathbb{T}}^{n}}\left( {\mathrm{{AF}}\left( {\mathbb{T}}^{n}\right) }\right)  \subseteq  \mathrm{{AF}}\left( {\mathbb{T}}^{n}\right)\), therefore $$\mathrm{{AF}}\left( {\mathbb{T}}^{n}\right)  = N\left( {\mathrm{{AF}}\left( {\mathbb{T}}^{n}\right) }\right)  = {E}_{{\mathbb{T}}^{n}}\left( {\mathrm{{AF}}\left( {\mathbb{T}}^{n}\right) }\right).$$ \hfill $\square$

\end{document}